\documentclass[12pt]{amsart}
\usepackage{amsmath,amssymb}
\newtheorem{theorem}{Theorem}
\newtheorem{lemma}{Lemma}
\newtheorem{proposition}{Proposition}
\begin{document}
\title{}
\author{}
\centerline{\Large Minimality in CR geometry and the CR Yamabe}
\vskip 0.1in \centerline{\Large problem on CR manifolds with
boundary}  \vskip 1cm  \centerline{\large Sorin
Dragomir\footnote{Author's address: Universit\`a degli Studi della
Basilicata, Dipartimento di Matematica, Campus Macchia Romana,
85100 Potenza, Italy, e-mail: {\tt dragomir@unibas.it}. The Author
acknowledges support from INdAM, Italy, within the
interdisciplinary project {\em Nonlinear subelliptic equations of
variational origin in contact geometry}.}}
\begin{abstract}
We study the minimality of an isometric immersion of a Riemannian
manifold into a strictly pseudoconvex CR manifold $M$ endowed with
the Webster metric (associated to a fixed contact form on $M$),
hence formulate a version of the CR Yamabe problem for CR
manifolds-with-boundary. This is shown to be a nonlinear
subelliptic problem of variational origin.
\end{abstract}
\maketitle

\section{Introduction}
Minimal surfaces $N^2$ in the lowest dimensional Heisenberg group
${\mathbb H}_1$, or more generally in a $3$-dimensional
nondegenerate CR manifold, have been recently considered by a
number of people (cf. N. Arcozzi \& F. Ferrari, \cite{kn:ArFe}, I.
Birindelli \& E. Lanconelli, \cite{kn:BiLa}, J-H. Cheng et alt.,
\cite{kn:CHMY}, N. Garofalo \& S.D. Pauls, \cite{kn:GaPa}, and
S.D. Pauls, \cite{kn:Pau}) motivated by the interest in a
Heisenberg version of the Bernstein problem, or by anticipating an
appropriate formulation of the CR Yamabe problem on a CR
manifold-with-boundary and a CR analog to the positive mass
theorem. All the notions of minimality dealt with are but ordinary
minimality of $N^2$ with respect to the ambient Webster metric.
This is demonstrated by our Theorem \ref{c:1} (though confined to
the case where the characteristic direction $T =
\partial /\partial t$ of ${\mathbb H}_1$ is tangent to $N^2$). We
also study minimality of a given isometric immersion $\Psi : N^m
\to {\mathbb H}_n$ of a $m$-dimensional Riemannian manifold $(N^m
, g)$ into $({\mathbb H}_n , g_{\theta_0} )$ (the Heisenberg group
carrying the Webster metric $g_{\theta_0}$ associated with the
contact form $\theta_0 = d t + i \sum_{j=1}^n (z_j d
\overline{z}^j - \overline{z}_j d z^j )$), cf. our Theorem
\ref{t:min}. A first step towards a Weierstrass type
representation of minimal surfaces in ${\mathbb H}_n$ is taken in
Theorem \ref{t:7}.
\par
The Yamabe problem on a compact $n$-dimensional ($n \geq 3$)
Riemannian manifold $(M , g)$ with boundary $\partial M$ is to
deform conformally the given metric $\hat{g} = u^{4/(n-2)} g$ ($u
> 0$) such that $(M , \hat{g})$ has constant scalar curvature and
$\partial M$ is minimal in $(M , \hat{g})$. This is equivalent to
solving the boundary value problem
\begin{equation}
\Delta u - \frac{n-2}{4(n-1)} \rho_g u + C u^{(n+2)/(n-2)} = 0
\;\; {\rm in} \;\; M, \label{e:1}
\end{equation}
\begin{equation}
\frac{\partial u}{\partial \eta} + \frac{n-2}{2} h_g u = 0 \;\;
{\rm on} \;\;
\partial M,
\label{e:2}
\end{equation}
where $\Delta$ and $\rho_g$ are respectively the Laplace-Beltrami
operator and the scalar curvature of $(M , g)$, $h_g$ is the mean
curvature of $\partial M \hookrightarrow (M , g)$, and $\eta$ is a
unit outward normal on $\partial M$ with respect to $g$. When $M$
is closed ($\partial M = \emptyset$) the full solution to
(\ref{e:1}) is described in \cite{kn:LePa}. When $\partial M \neq
\emptyset$ the problem (\ref{e:1})-(\ref{e:2}) was solved by J.F.
Escobar, \cite{kn:Esc}, under the assumptions that 1) $n \in \{
3,4,5\}$, or 2) $n \geq 3$ and $\partial M$ has some nonumbilic
point, or 3) $n \geq 6$, $\partial M$ is totally umbilical, and
either $M$ is locally conformally flat or the Weyl tensor doesn't
vanish identically on $\partial M$. A CR analog of the Yamabe
problem was formulated by D. Jerison \& J.M. Lee, \cite{kn:JL1},
though only on closed CR manifolds. Precisely, if $M$ is a
$(2n+1)$-dimensional closed strictly pseudoconvex CR manifold on
which a contact form $\theta$ has been fixed then the CR Yamabe
problem is to look for a contact form $\hat{\theta} = u^{p-2}
\theta$ ($p = 2 + 2/n$) such that the Tanaka-Webster connection of
$(M , \hat{\theta})$ has constant pseudohermitian scalar curvature
$\hat{\rho} = \lambda$. This is equivalent to solving
\begin{equation}
- (2 + 2/n) \, \Delta_b u + \rho \, u = \lambda \, u^{p-1}
\label{e:3}
\end{equation}
(the {\em CR Yamabe equation}) where $\Delta_b$ and $\rho$ are
respectively the sublaplacian\footnote{As to the sign convention
the sublaplacian in \cite{kn:JL2} is $-\Delta_b$.} and the
pseudohermitian scalar curvature of $(M , \theta )$. D. Jerison \&
J.M. Lee solved (cf. \cite{kn:JL2}-\cite{kn:JL3}) the problem
(\ref{e:3}) under the assumption that\footnote{If $n \geq 2$ and
$M$ is not locally CR equivalent to $S^{2n+1}$ then $\lambda (M) <
\lambda (S^{2n+1})$, cf. \cite{kn:JL2}.} $\lambda (M) < \lambda
(S^{2n+1})$, where $\lambda (M)$ is the CR invariant
\[ \inf \{ \int_M ( b_n \| \pi_H \nabla u \|^2 + \rho u^2 ) \theta
\wedge (d \theta )^n : \int_M |u|^p \theta \wedge (d \theta )^n =
1 \} . \] Moreover, the inequality $\lambda (M) \leq \lambda
(S^{2n+1})$ holds true. The remaining case $\lambda (M) = \lambda
(S^{2n+1})$ was settled by N. Gamara \& R. Yacoub, \cite{kn:GaYa}.
It is noteworthy that the proof in \cite{kn:GaYa} doesn't rely on
a CR analog to the positive mass theorem, but rather on techniques
within the theory of critical points at infinity (by analogy with
A. Bahri \& H. Brezis, \cite{kn:BaBr}). When $\partial M \neq
\emptyset$ no formulation of the CR Yamabe problem is available as
yet, perhaps due to the previous lack of a natural CR analog to
minimality.
\par
Our approach (as well as in \cite{kn:JL2}) is to formulate the CR
Yamabe problem as the Yamabe problem for the Fefferman metric
$F_\theta$, a Lorentz metric on the total space $C(M)$ of the
canonical circle bundle $S^1 \to C(M) \stackrel{\pi}{\rightarrow}
M$ (cf. \cite{kn:Lee}). That is, to look for a positive function
$u \in C^\infty (M)$ such that the Fefferman metric
$F_{\hat{\theta}}$ corresponding to the contact form $\hat{\theta}
= u^{p-2} \theta$ has constant scalar curvature. What is the
appropriate boundary condition?
\par
When $\partial M$ is nonempty $C(M)$ is a manifold-with-boundary
as well, and (by Theorem \ref{p:5}) the tangent space $T_z
(\partial C(M))$ is nondegenerate in $(T_z (C(M)), F_{\theta ,
z})$ at all points $z$, except for those projecting on ${\rm
Sing}(T^T )$, the singular points of the tangential component
(with respect to $\partial M$) of the characteristic direction $T$
of $d \theta$. It also turns out that $\partial C(M) \setminus
\pi^{-1}({\rm Sing}(T^T ))$ is a Lorentz manifold (with the metric
induced by $F_\theta$). Therefore, when ${\rm Sing}(T^T ) =
\emptyset$ we may request that $\partial C(M)$ be minimal in
$(C(M), F_{\hat{\theta}})$. By Theorem \ref{t:1} this projects to
the natural boundary condition (\ref{e:34}) on $\partial M$, thus
leading to the {\em CR Yamabe problem} (\ref{e:33})-(\ref{e:34})
on a CR manifold-with-boundary. This is shown (cf. Theorem
\ref{p:11}) to be a nonlinear subelliptic problem of variational
origin. \vskip 0.1in {\small {\bf Acknowledgements}. The Author is
grateful to E. Lanconelli for stimulating conversations on the
arguments in this paper and for introducing him to the results in
the preprint \cite{kn:BiLa}. Also, the Author wishes to express
his gratitude for the hospitality and excellent working atmosphere
in the Department of Mathematics of the University of Bologna and
for discussions with N. Arcozzi and F. Ferrari (who kindly
provided the preprint \cite{kn:ArFe}).}

\section{CR manifolds with boundary}
Let $M$ be an oriented $m$-dimensional $C^\infty$
manifold-with-boundary $\partial M$. A {\em CR structure} is a
complex subbundle $T_{1,0}(M)$ of the complexified tangent bundle
$T(M) \otimes {\mathbb C}$, of complex rank $n$ ($0 < n \leq
[m/2]$), such that
\[ T_{1,0}(M) \cap T_{0,1}(M) = (0), \]
\[ Z,W \in \Gamma^\infty (T_{1,0}(M)) \Longrightarrow [Z,W] \in
\Gamma^\infty (T_{1,0}(M)). \] Here $T_{0,1}(M) =
\overline{T_{1,0}(M)}$ (complex conjugation). The pair $(M ,
T_{1,0}(M))$ is a {\em CR manifold} (with boundary) and the
integer $n$ is its {\em CR dimension}. Also $k = m - 2n$ is its
{\em CR codimension} and the pair $(n,k)$ is its {\em type}.
\par
There is a natural first order differential operator
$\overline{\partial}_b$ (the {\em tangential Cauchy-Riemann
operator}) given by $(\overline{\partial}_b u) \overline{Z} =
\overline{Z}(u)$, for any $C^1$ function $u : M \to {\mathbb C}$
and any $Z \in T_{1,0}(M)$. Then $\overline{\partial}_b u = 0$ are
the {\em tangential Cauchy-Riemann equations}. A solution to the
tangential Cauchy-Riemann equations is a {\em CR function} on $M$.
Let ${\rm CR}^r (M)$ denote the space of all CR functions on $M$
of class $C^r$.
\par
The boundary $\partial M$ is {\em noncharacteristic} for
$T_{1,0}(M)$ if for any local frame $\{ T_\alpha : 1 \leq \alpha
\leq n \}$ of $T_{1,0}(M)$ defined on the open subset $U \subseteq
M$ one has $T_\alpha \not\in T(\partial M) \otimes {\mathbb C}$
(i.e. $T_{\alpha , x} \not\in T_x (\partial M) \otimes {\mathbb
C}$, for some $x \in U \cap
\partial M$) for some $1 \leq \alpha \leq n$.
\par
The {\em Levi distribution} of the CR manifold $(M , T_{1,0}(M))$
is
\[ H(M) = {\rm Re}\{ T_{1,0}(M) \oplus T_{0,1}(M) \} . \]
It carries the complex structure \[ J : H(M) \to H(M), \;\; J(Z +
\overline{Z}) = i (Z - \overline{Z}), \;\; Z \in T_{1,0}(M). \]
Assume from now on that $M$ is a CR manifold of type $(n,1)$ (of
{\em hypersurface type}). $H(M)$ is oriented by $J$, hence the
conormal bundle
\[ H(M)^\bot_x = \{ \omega \in T^*_x (M) : {\rm Ker}(\omega )
\supseteq H(M)_x \} , \;\;\; x \in M, \] is an oriented real line
bundle, hence trivial. Let then $\theta$ be a global nowhere
vanishing section in $H(M)^\bot$ (a {\em pseudohermitian
structure} on $M$). The {\em Levi form} is
\[ L_\theta (Z, \overline{W}) = - i (d \theta )(Z, \overline{W}),
\;\;\; Z, W \in T_{1,0}(M), \] and $M$ is {\em nondegenerate}
(respectively {\em strictly pseudoconvex}) if $L_\theta$ is
nondegenerate (respectively positive definite) for some $\theta$.
Also $M$ is {\em Levi flat} if $L_\theta = 0$ (equivalently, if
$H(M)$ is integrable). An alternative definition of the Levi form
is \[ G_\theta (X,Y) = (d \theta )(X, J Y), \;\;\; X,Y \in H(M).
\] Note that $L_\theta$ and the ${\mathbb C}$-linear extension of $G_\theta$
coincide on $T_{1,0}(M) \otimes T_{0,1}(M)$. If $M$ is
nondegenerate then any pseudohermitian structure $\theta$ is a
{\em contact form}, i.e. $\theta \wedge (d \theta )^n$ is a volume
form on $M$. Let $M$ be a nondegenerate CR manifold and $\theta$ a
fixed contact form (the pair $(M, \theta )$ is commonly referred
to as a {\em pseudohermitian manifold}). There is a unique vector
field $T$ on $M$ such that $\theta (T) = 1$ and $(d \theta )(T ,
X) = 0$, for any $X \in T(M)$ ($T$ is the {\em characteristic
direction} of $d \theta$). The {\em Webster metric} of $(M ,
\theta )$ is given by
\[ g_\theta (X,Y) = G_\theta (X , J Y), \;\; g_\theta (X,T) = 0,
\;\; g_\theta (T,T) = 1, \] for any $X,Y \in H(M)$. $g_\theta$ is
a semi-Riemannian (Riemannian, if $M$ is strictly pseudoconvex and
$L_\theta$ is positive definite) metric on $M$.
\begin{proposition} Let $M$ be a nondegenerate CR manifold-with-boundary. Then the
boundary $\partial M$ is noncharacteristic for $T_{1,0}(M)$.
\label{p:1}
\end{proposition}
\noindent The proof is by contradiction. Assume that there is a
local frame $\{ T_\alpha \}$ of $T_{1,0}(M)$ on $U \subseteq M$
such that $T_\alpha \in T(\partial M) \otimes {\mathbb C}$, for
all $1 \leq \alpha \leq n$. Then $T_{1,0}(M)_x \subset T_x
(\partial M) \otimes {\mathbb C}$ for any $x \in U \cap \partial
M$. Then, by taking complex conjugates, $T_{0,1}(M)_x \subset T_x
(\partial M) \otimes {\mathbb C}$ hence, by looking at dimensions,
$H(M)_x = T_x (\partial M)$, i.e. $L_{\theta , x} = 0$, a
contradiction. $\square$ \vskip 0.1in From now on we assume that
$M$ is nondegenerate. For each boundary point $x \in
\partial M$ we set
\[ T_{1,0}(\partial M)_x = T_{1,0}(M)_x \cap [T_x (\partial M)
\otimes {\mathbb C}]. \] Let $\{ T_\alpha : 1 \leq \alpha \leq n
\}$ be a local frame of $T_{1,0}(M)$, defined on the local
coordinate neighborhood $(U , \varphi = (x^1 , \cdots ,
x^{2n+1}))$. $U \cap \partial M$ consists of the points $x \in U$
such that $\varphi (x) \in \partial {\mathbb R}^{2n+1}_+ =
{\mathbb R}^{2n} \times \{ 0 \}$. We may write $T_\alpha =
f^A_\alpha \; \partial /\partial x^A$,  for some $C^\infty$
functions $f^A_\alpha : U \to {\mathbb C}$. By Proposition
\ref{p:1} there is $\alpha$, say $\alpha = n$, such that $T_\alpha
\not\in T(\partial M) \otimes {\mathbb C}$. Then $f_n^{2n+1}(x_0 )
\neq 0$ for some $x_0 \in U \cap \partial M$, and then $f^{2n+1}_n
\neq 0$ on a whole neighborhood of $x_0$, which we may denote
again by $U$. Then
\[ \{ T_j - \left( \lambda^{2n+1}_j /\lambda^{2n+1}_n \right) \; T_n : 1
\leq j \leq n-1 \} \] is a local frame of $T_{1,0}(\partial M)$ on
$U \cap \partial M$, hence $T_{1,0}(\partial M)$ has rank $n-1$.
We got
\begin{proposition} Let $M$ be a nondegenerate CR manifold-with-boundary,
of CR dimension $n$. Then its boundary $\partial M$ is a
CR manifold of type $(n-1,2)$, i.e. $T_{1,0}(\partial M) =
T_{1,0}(M) \cap [T(\partial M) \otimes {\mathbb C}]$ is a CR
structure of CR codimension $2$. \label{p:2}
\end{proposition}
\noindent Let us look at a few examples. For instance, let
${\mathbb H}_n = {\mathbb C}^n \times {\mathbb R}$ be the
Heisenberg group, with the CR structure spanned by
\[ Z_j = \frac{\partial}{\partial z_j} + i \, \overline{z}_j
\, \frac{\partial}{\partial t}\, , \;\;\; 1 \leq j \leq n, \] (if
$n=1$ then $\overline{Z}_1$ is the {\em Lewy operator}, cf.
\cite{kn:Lew}). ${\mathbb H}_n$ is a Lie group with the group law
\[ (z, t) \cdot (w , s) = (z + w , t + s + 2 \; {\rm Im}(z \cdot
\overline{w})), \] for $(z,t), \; (w, t) \in {\mathbb H}_n$, where
$z \cdot \overline{w} = \delta_{jk} z^j \overline{w}^k$ (with the
convention $z^j = z_j$), and $Z_j$ are left invariant. \vskip
0.1in\noindent {\bf Example 1.} ${\mathbb H}_n^+ = \{ (z,t) \in
{\mathbb H}_n : t \geq 0 \}$ is a CR manifold-with-boundary
$\partial {\mathbb H}_n^+ = {\mathbb C}^n \times \{ 0 \}$. Let $U
= \{ (z,t) \in {\mathbb H}_n^+ : z_n \neq 0 \}$. Then
\begin{equation}
\{ \frac{\partial}{\partial z^a} -
\frac{\overline{z}_a}{\overline{z}_n} \; \frac{\partial}{\partial
z_n} : 1 \leq a \leq n-1 \}
\label{e:4}
\end{equation}
is a local frame of $T_{1,0}(\partial {\mathbb H}_n^+ )$ on $U
\cap \partial {\mathbb H}_n^+$. In particular, the tangential
Cauchy-Riemann equations on $\partial {\mathbb H}^+_n$ are
\[
z_n \; \frac{\partial u}{\partial \overline{z}_a}  - z_a \;
\frac{\partial u}{\partial \overline{z}_n} = 0, \;\;\; 1 \leq a
\leq n-1.
\]
$\square$ \vskip 0.1in\noindent The {\em Heisenberg norm} is $|x|
= (|z|^4 + t^2 )^{1/4}$, for any $x = (z,t) \in {\mathbb H}_n$,
where $|z|^2 = z \cdot \overline{z}$. \vskip 0.1in\noindent {\bf
Example 2.} $\Omega_r = \{ x \in {\mathbb H}_n : |x| \leq r \}$
($r
> 0$) is a CR manifold-with-boundary $\partial \Omega_r = \Sigma_r
= \{ x \in {\mathbb H}_n : |x| = r \}$ (the {\em Heisenberg
sphere}, cf. \cite{kn:GaLa}). Let us set $\phi (z,t) = |z|^2 - i
\; t$. Note that $\overline{\partial}_b \phi = 0$, i.e. $\phi \in
{\rm CR}^\infty ({\mathbb H}_n )$. Taking into account that
\[ Z_j (|x|) = \frac{\overline{\phi}}{2\, |x|^3} \;
\overline{z}_j \, , \;\;\; 1 \leq j \leq n, \] it follows that
(\ref{e:4}) is a local frame of $T_{1,0}(\Sigma_r )$ on $\Sigma_r
\cap \{ z \in {\mathbb H}_n : z_n \neq 0 \}$. The {\em
Folland-Stein operators} are
\begin{equation}
{\mathcal L}_\alpha = - \frac{1}{2} \sum_{j=1}^n (Z_j
\overline{Z}_j + \overline{Z}_j Z_j ) + i \; \alpha \; T, \;\;\;
\alpha \in {\mathbb C}, \label{e:folste}
\end{equation}
where $T = \partial /\partial t$. Let us consider the function
\[ \varphi_\alpha (z,t) = \phi (z,t)^{- (n+\alpha )/2} \;
\overline{\phi (z,t)}^{\; - (n-\alpha )/2} \, , \] and the
constant $c_\alpha = 2^{2-2n} \pi^{n+1}/\left( \Gamma
(\frac{n+\alpha}{2}) \; \Gamma (\frac{n-\alpha}{2})\right)$.
$\alpha \in {\mathbb C}$ is {\em admissible} if $c_\alpha \neq 0$
(equivalently if $\pm \alpha \in \{ n , \, n+2, \, n+4 , \, \cdots
\}$). The Folland-Stein operators (\ref{e:folste}) form a family
of operators of the form $A + \alpha B$ (where $A$ is a second
order hypoelliptic operator and $B$ is a first order operator)
which are hypoelliptic for any admissible $\alpha$ (cf.
\cite{kn:FoSt}, p. 444). This is by now classical, and as well
known the key ingredient in the proof is to build a fundamental
solution to (\ref{e:folste}) i.e. to show that ${\mathcal
L}_\alpha ( \varphi_\alpha /c_\alpha ) = \delta$, for any
admissible $\alpha$. It is noteworthy that the Heisenberg spheres
$\Sigma_r$ are the level sets of
\[ \varphi_0 (z,t) = |\phi (z,t)|^{-n} = \left( |z|^4 + t^2
\right)^{-n/2}\,  .  \] Let $\theta_0$ be the canonical
pseudohermitian structure on ${\mathbb H}_n$ i.e.
\[ \theta_0 = d t + i \sum_{j=1}^n (z_j d \overline{z}^j -
\overline{z}_j d z^j ). \] ${\mathbb H}_n$ is strictly
pseudoconvex and $L_{\theta_0}$ is positive definite. Moreover,
the Webster metric of $({\mathbb H}_n , \theta_0 )$ is expressed
by
\[ g_{\theta_0}(X_j , X_k ) = g_{\theta_0}(Y_j , Y_k ) =
\delta_{jk} \, , \;\;\; g_{\theta_0}(X_j , Y_k ) = 0, \]
\[ g_{\theta_0}(X_j , T) = g_{\theta_0} (Y_j , T) = 0, \;\;\;
g_{\theta_0}(T,T) = 1, \] where
\[ X_j = \frac{1}{\sqrt{2}} \, (Z_j + \overline{Z}_j ), \;\;\; Y_j
= \frac{i}{\sqrt{2}} \, (Z_j - \overline{Z}_j ). \]
\begin{proposition} The Heisenberg spheres form a foliation of
$({\mathbb H}_n , g_{\theta_0})$ whose normal bundle is the span
of \begin{equation} V = T + (\phi /t) z^j Z_j +
(\overline{\phi}/t) \overline{z}^j \overline{Z}_j \, .
\label{e:6}
\end{equation} \label{p:fol}
\end{proposition}
\noindent Then perhaps (\ref{e:6}) is the Heisenberg analog to the
radial vector field in ${\mathbb R}^{2n+1}$ (see \cite{kn:GaLa},
p. 331-332). \par\noindent {\em Proof of Proposition} \ref{p:fol}.
Let us set
\[ E_j = Z_j + \overline{Z}_j - \frac{1}{t} (\phi z_j +
\overline{\phi} \overline{z}_j ) T, \;\; F_j = i(Z_j -
\overline{Z}_j ) + \frac{i}{t} (\phi z_j - \overline{\phi}
\overline{z}_j ) T. \] Then $\{ E_j , F_j \}$ is a local frame of
the tangent bundle of the foliation and a calculation shows that
(\ref{e:6}) satisfies $g_{\theta_0} (E_j , V) = g_{\theta_0}(F_j ,
V) = 0$.  $\square$ \vskip 0.1in Let $M$ and $N$ be two CR
manifolds with boundary. A {\em CR map} is a $C^\infty$ map $f : M
\to N$ such that $(d_x f) T_{1,0}(M)_x \subseteq
T_{1,0}(N)_{f(x)}$, for any $x \in M$. A {\em CR immersion} is an
immersion and a CR map. A CR immersion $f : M \to N$ is {\em neat}
if i) $f(M) \cap
\partial N = f(\partial M)$ and ii) for each point $x \in \partial
M$ there is a local chart $\psi : V \to {\mathbb R}^{m+p}_+$ of
$N$ such that $f(x) \in V$ and $\psi^{-1}({\mathbb R}^m_+ ) = V
\cap f(M)$ ($m = \dim (M)$). \vskip 0.1in\noindent {\bf Example
3.} $\Sigma_r^+ = \Sigma_r \cap {\mathbb H}^+_n$ is a CR
manifold-with-boundary $\partial \Sigma_r^+ = S^{2n-1}(r) \times
\{ 0 \}$ and the inclusion $\Sigma^+_r \to {\mathbb H}^+_n$ is a
neat CR immersion. $\square$ \vskip 0.1in\noindent {\bf Example
4.} $S^{2n+1}_+ = S^{2n+1} \cap {\mathbb R}^{2n+2}_+$ is a CR
manifold-with-boundary $\partial S^{2n+1}_+ = S^{2n} \times \{ 0
\}$. Let ${\mathcal C}$ be the Cayely transform
\[ {\mathcal C}(\zeta ) = (\frac{\zeta^\prime}{1 + \zeta^{n+1}} \, , \, i \, \frac{1 -
\zeta^{n+1}}{1 + \zeta^{n+1}}), \;\; \zeta = (\zeta^\prime ,
\zeta^{n+1}), \;\; 1 + \zeta^{n+1} \neq 0, \] and $f : {\mathbb
H}_n \to \partial \Omega_{n+1}$ the CR isomorphism $f(z,t) = (z,
t+ i |z|^2 )$ with the obvious inverse $f^{-1}(z,w) = (z, {\rm
Re}(w))$. Here $\Omega_{n+1}$ is the Siegel domain
\[ \Omega_{n+1} = \{ (z, w) \in {\mathbb C}^{n+1} : {\rm Im}(w) >
|z|^2 \} . \] Then $F = f^{-1} \circ {\mathcal C}$ is a neat CR
diffeomorphism  \[ F : S^{2n+1}_+ \setminus \{ (0, \cdots , 0,
-1)\} \to {\mathbb H}_n^+ . \] Indeed if $\zeta \in S^{2n+1}_+$
and $\zeta^{n+1} = u + iv$ ($v \geq 0$) and $(z,t) = F(\zeta )$
then $t = 2v/[(1+u)^2 + v^2 ] \geq 0$. In particular $F$ descends
to a CR diffeomorphism $(S^{2n} \times \{ 0 \}) \setminus \{ (0,
\cdots , 0 , -1) \} \approx {\mathbb C}^n \times \{ 0 \}$.
$\square$ \vskip 0.1in Let $M$ be a nondegenerate CR
manifold-with-boundary. A complex $p$-form $\eta$ on $M$ is a
$(p,0)$-{\em form} if $T_{0,1}(M) \, \rfloor \, \eta = 0$. Let
$\Lambda^{p,0}(M) \to M$ be the bundle of all $(p,0)$-forms. If
$M$ has CR dimension $n$ then the top degree $(p,0)$-forms are the
$(n+1,0)$-forms. $K(M) = \Lambda^{n+1,0}(M)$ is the {\em canonical
bundle} over $M$. There is a natural action of ${\mathbb R}_+ =
(0, + \infty )$ on $K(M) \setminus \{ 0 \}$. Let $C(M)$ be the
quotient space and $\pi : C(M) \to M$ the projection. This
construction leads to a principal bundle $S^1 \to C(M) \to M$ (the
{\em canonical circle bundle} over $M$). Let $\theta$ be a
pseudohermitian structure on $M$ and $T$ the characteristic
direction of $d \theta$. Given a local frame $\{ T_\alpha \}$ of
$T_{1,0}(M)$ on a local coordinate neighborhood $(U, x^A )$, let
$\theta^\alpha$ be the locally defined complex $1$-forms
determined by
\[ \theta^\alpha (T_\beta ) = \delta_\beta^\alpha \, , \;\;
\theta^\alpha (T_{\overline{\beta}} ) = 0, \;\; \theta^\alpha (T)
= 0. \] Here $T_{\overline{\alpha}} = \overline{T}_\alpha$. Then
\[ \pi^{-1}(U) \to U \times S^1 \, , \;\; [z] \mapsto
(x \, , \, \lambda /|\lambda | ), \]
\[ z = \lambda \, (\theta \wedge \theta^1 \wedge \cdots \wedge
\theta^n )_x \, , \;\; \lambda \in {\mathbb C} \setminus \{ 0 \} ,
\;\; x \in M, \] is a local trivialization chart of the canonical
circle bundle. Let us set $\gamma : \pi^{-1} (U) \to {\mathbb R}$,
$\gamma ([z]) = \arg (\lambda )$ (where $\arg : {\mathbb C} \to
[0, 2 \pi )$). Then $(\pi^{-1}(U) , \, \tilde{x}^A = x^A \circ \pi
, \, \gamma )$ are naturally induced local coordinates on $C(M)$
and $\pi^{-1}(U \cap \partial M)$ consists of all $c \in
\pi^{-1}(U)$ with $\tilde{x}^{2n+1}(c) = 0$, i.e. $C(M)$ is a
manifold-with-boundary modelled on ${\mathbb R}^{2n+1}_+ \times
{\mathbb R}$. We obtained
\begin{lemma}
Let $M$ be a nondegenerate CR manifold-with-boundary. Then the
total space $C(M)$ of the canonical circle bundle is a
manifold-with-boundary $\partial C(M) = \pi^{-1}(\partial M)$. In
particular $\partial C(M)$ is a principal $S^1$-bundle over
$\partial M$.
\end{lemma}
\noindent Let $\nabla$ be the unique linear connection on $M$ (the
{\em Tanaka-Webster connection}) satisfying the axioms 1) $H(M)$
is parallel with respect to $\nabla$, 2) $\nabla J = 0$, $\nabla
g_\theta = 0$, and 3) the torsion $T_\nabla$ of $\nabla$ is {\em
pure}, i.e. $T_\nabla (Z, W) = 0$, $T_\nabla (Z, \overline{W}) = 2
i G_\theta (Z , \overline{W}) T$, and $\tau \circ J + J \circ \tau
= 0$. Here $\tau (X) = T_\nabla (T ,X)$ is the {\em
pseudohermitian torsion}. We set $A(X,Y) = g_\theta (\tau X, Y)$,
for any $X,Y \in T(M)$. By a result of S. Webster, \cite{kn:Web},
$A$ is symmetric.
\par
With respect to a local frame $\{ T_\alpha : 1 \leq \alpha \leq n
\}$ of $T_{1,0}(M)$, defined on an open set $U \subseteq M$, it is
customary to set $g_{\alpha\overline{\beta}} = L_\theta (T_\alpha
, T_{\overline{\beta}})$ (the local coefficients of the Levi
form), $\nabla T_\beta = \omega_\beta^\alpha \otimes T_\alpha$
(the connection $1$-forms) and $R^\nabla (T_A , T_B ) T_C =
{{R_C}^D}_{AB} T_D$ (the curvature components). The range of the
indices $A,B,C, \cdots$ is $\{ 0, 1, \cdots , n , \overline{1},
\cdots , \overline{n} \}$ (with the convention $T_0 = T$). Next,
the {\em pseudohermitian Ricci tensor} is
$R_{\lambda\overline{\mu}} =
{{R_\lambda}^\alpha}_{\alpha\overline{\mu}}$ and the {\em
pseudohermitian scalar curvature} is $\rho =
g^{\alpha\overline{\beta}} R_{\alpha\overline{\beta}}$. When $M$
is strictly pseudoconvex and $\theta$ is a pseudohermitian
structure such that $L_\theta$ is positive definite $C(M)$ carries
a Lorentz metric $F_\theta$ such that $F_{\hat{\theta}} = e^{u
\circ \pi} F_\theta$, where $\hat{\theta} = e^u \theta$, $u \in
C^\infty (M)$ (in particular the {\em restricted conformal class}
$[F_\theta ] = \{ e^{u \circ \pi} F_\theta : u \in C^\infty (M)
\}$ is a CR invariant). Cf. J.M. Lee, \cite{kn:Lee}, $F_\theta$ is
given by
\begin{equation} F_\theta =
\pi^* \tilde{G}_\theta + 2 (\pi^* \theta ) \odot \sigma ,
\label{e:ad2}
\end{equation}
\begin{equation} \label{e:ad1} \sigma = \frac{1}{n+2} \{ d \gamma + \pi^* (i \,
\omega^\alpha_\alpha - \frac{i}{2} \, g^{\alpha\overline{\beta}} d
g_{\alpha\overline{\beta}} - \frac{\rho}{4(n+1)}  \, \theta ) \} .
\end{equation}
$F_\theta$ is the {\em Fefferman metric} of $(M , \theta )$. Here
$\tilde{G}_\theta$ is the degenerate $(0,2)$-tensor field on $M$
given by
\[ \tilde{G}_\theta (X,Y) = (d \theta )(X, J Y), \;\;\;
\tilde{G}_\theta (T , Z) = 0, \] for any $X, Y \in H(M)$ and any
$Z \in T(M)$. Also $\odot$ denotes the symmetric tensor product.
\par
Let $S = \partial /\partial \gamma$ be the tangent to the
$S^1$-action. $\sigma$ is a connection $1$-form in $S^1 \to C(M)
\to M$. If $X \in T(M)$ is a tangent vector field on $M$ then
$X^\uparrow \in T(C(M))$ will denote the horizontal lift of $X$
with respect to the connection ${\mathcal H} = {\rm Ker}(\sigma
)$. Although the submersion $\pi : C(M) \to M$ is not
semi-Riemannian (its fibres are degenerate) a technique similar to
that in \cite{kn:One} leads to
\begin{lemma} For any $X,Y \in H(M)$
\[ \nabla^{C(M)}_{X^\uparrow} Y^\uparrow = (\nabla_X Y)^\uparrow -
(d \theta )(X,Y) T^\uparrow - (A(X,Y) + (d \sigma )(X^\uparrow ,
Y^\uparrow )) \hat{S}, \]
\[ \nabla^{C(M)}_{X^\uparrow} T^\uparrow  = (\tau X + \phi
X)^\uparrow , \]
\[ \nabla^{C(M)}_{T^\uparrow} X^\uparrow = (\nabla_T X + \phi
X)^\uparrow + 2(d \sigma )(X^\uparrow , T^\uparrow ) \hat{S},  \]
\[ \nabla^{C(M)}_{X^\uparrow} \hat{S} = \nabla^{C(M)}_{\hat{S}} X^\uparrow = (J
X)^\uparrow , \]
\[ \nabla^{C(M)}_{T^\uparrow} T^\uparrow = V^\uparrow , \;\;
\nabla^{C(M)}_{\hat{S}} \hat{S} = 0, \]
\[ \nabla^{C(M)}_{\hat{S}} T^\uparrow = \nabla^{C(M)}_{T^\uparrow} \hat{S} = 0,
\]
where $\phi : H(M) \to H(M)$ is given by $G_\theta (\phi X , Y) =
(d \sigma )(X^\uparrow , Y^\uparrow )$, and $V \in H(M)$ is given
by $G_\theta (V , Y) = 2 (d \sigma )(T^\uparrow , Y^\uparrow )$.
Also $\hat{S} = ((n+2)/2) S$. \label{l:1}
\end{lemma}
Lemma \ref{l:1} relates the Levi-Civita connection $\nabla^{C(M)}$
of $(C(M), F_\theta )$ to the Tanaka-Webster connection of $(M,
\theta )$. Cf. \cite{kn:BaDrUr} for the proof of Lemma \ref{l:1}.

\section{The geometry of the first fundamental form of the
boundaries} Let $M$ be a strictly pseudoconvex CR manifold and
$\theta$ a contact form on $M$ such that $G_\theta$ is positive
definite. Let $T(\partial M)^\bot \to \partial M$ be the normal
bundle of $\partial M \hookrightarrow (M , g_\theta )$. Let ${\rm
tan}_x : T_x (M) \to T_x (\partial M)$ and ${\rm nor}_x : T_x (M)
\to T(\partial M)^\bot_x$ be the projections associated with the
direct sum decomposition
\[ T_x (M) = T_x (\partial M) \oplus T(\partial M)^\bot_x \, ,
\;\;\; x \in \partial M. \] If $T$ is the characteristic direction
of $d \theta$ then we set $T^\bot = {\rm nor} (T)$ and $T^T = {\rm
tan}(T)$.
\begin{theorem}  Let ${\rm Null}(j^* F_\theta )$ consist of all
$V \in T(\partial C(M))$ such that $F_\theta (V,W) = 0$, for any
$W \in T(\partial C(M))$. Let us consider the closed set ${\rm
Sing}(T^T ) = \{ x \in \partial M : T^T_x = 0 \}$ and set $\Omega
=
\partial M \setminus {\rm Sing}(T^T )$. Then
\[ {\rm Null}(j^* F_\theta )_z = \begin{cases} 0, & z \in
\pi^{-1}(\Omega ), \cr {\rm Ker}(d_z \pi ), & z \in \pi^{-1}({\rm
Sing}(T^T )), \cr
\end{cases} \]
for any $z \in \partial C(M)$. Moreover $(\pi^{-1}(\Omega )\, , \,
j^* F_\theta )$ is a Lorentz manifold. \label{p:5}
\end{theorem}
\noindent Here $j : \partial C(M) \hookrightarrow C(M)$ is the
inclusion. Hence $\partial C(M)$ is degenerate at each point $z
\in \pi^{-1}({\rm Sing}(T^T ))$. In particular, if $\partial M$ is
tangent to $T$ then the boundary $(\partial C(M) \, , \, j^*
F_\theta )$ is a Lorentz manifold. \vskip 0.1in\noindent {\bf
Example 1.} ({\em continued})
\par\noindent  $T(\partial {\mathbb H}_n^+ )$ is the span of
$\{ X_j - \sqrt{2}\,  y_j T , \; Y_j + \sqrt{2} \, x_j T : 1 \leq
j \leq n \}$ hence $\xi = T + \sqrt{2} \, y^j X_j - \sqrt{2} \,
x^j Y_j$ is normal to $\partial {\mathbb H}_n^+$ (with $z^j = x^j
+ i y^j$). Then $T$ decomposes as
\[ T = a^j (X_j - \sqrt{2} \, y_j T) + b^j (Y_j + \sqrt{2} \, x_j T) + c
\xi , \]
\[ a^j = - \frac{\sqrt{2}\,  y^j}{1 + 2 |z|^2} \, , \;\; b^j =
\frac{\sqrt{2} \, x^j}{1 + 2 |z|^2} \, , \;\; c = \frac{1}{1 + 2
|z|^2} \, . \] Then $T^\bot = c \xi$ and (with the conventions in
Theorem \ref{p:5}) ${\rm Sing}(T^T ) = \{ 0 \}$ hence $(\partial
C({\mathbb H}_n^+ ) \setminus \pi^{-1} (0) \, , \,  j^*
F_{\theta_0})$ is a Lorentz manifold. $\square$ \vskip 0.1in {\em
Proof of Theorem} \ref{p:5}. Let $V \in T(\partial C(M))$ such
that $F_\theta (V , W) = 0$ for any $W \in T(\partial C(M))$ i.e.
\[ (\pi^* \tilde{G}_\theta )(V , W) + (\pi^* \theta )(V) \sigma
(W) + (\pi^* \theta )(W) \sigma (V) = 0. \] By taking into account
\begin{equation}
T(C(M)) = {\rm Ker}(\sigma ) \oplus {\rm Ker}(d \pi )
\label{e:old1}
\end{equation}
we may decompose $V = V_H + V_V$, with $V_H \in {\rm Ker}(\sigma
)$. Then
\begin{equation}
\tilde{G}((d \pi ) V_H , (d \pi )W_H ) + \theta ((d \pi ) V_H )
\sigma (W_V ) + \theta ((d \pi ) W_H ) \sigma (V_V ) = 0.
\label{e:old2}
\end{equation}
As $\partial C(M)$ is a saturated set, it is tangent to the
$S^1$-action. Hence we may apply (\ref{e:old2}) for $W = S \in
{\rm Ker}(d \pi ) \subset T(\partial C(M))$. As $\sigma (S) =
1/(n+2)$ we obtain
\[ \theta ((d \pi ) V_H ) = 0, \]
i.e. $(d \pi )V_H \in H(M)$, hence (\ref{e:old2}) becomes
\begin{equation}
\tilde{G}_\theta ((d \pi ) V_H , (d \pi )W_H ) + \theta ((d \pi )
W_H ) \sigma (V_V ) = 0. \label{e:old3}
\end{equation}
Applying (\ref{e:old3}) for $W = V$ gives
\[ G_\theta ((d \pi )V_H , (d \pi )V_H ) = 0 \]
hence $(d \pi ) V_H = 0$, and then $V_H = 0$ (due to ${\rm
Ker}(\sigma ) \cap {\rm Ker}(d \pi ) = (0)$). Therefore, on one
hand \begin{equation}
\label{e:incl}
{\rm Null}(j^* F_\theta )
\subseteq {\rm Ker}(d \pi )
\end{equation}
and on the other (\ref{e:old3}) becomes
\begin{equation}
\theta ((d \pi ) W_H ) \sigma (V_V ) = 0. \label{e:old4}
\end{equation}
Let $x_0 \in \Omega$ (so that $T^T_{x_0} \neq 0$) and $z_0 \in
\pi^{-1}(x_0 )$. We may apply (\ref{e:old4}) for $W = (T^T
)^\uparrow$, at the point $z_0$. Yet
\[ (\pi^* \theta )(W_H )_{z_0} = \theta (T^T )_{x_0} = \| T^T
\|^2_{x_0} \neq 0 \] hence (by (\ref{e:old4})) $\sigma (V_V
)_{z_0} = 0$, or $(V_V )_{z_0} = 0$ and we may conclude that ${\rm
Null}(j^* F_\theta )_{z_0} = (0)$. To complete the proof of
Theorem \ref{p:5} it suffices to show that ${\rm Null}(j^*
F_\theta )_{z}$ is $1$-dimensional, for any $z \in \pi^{-1} (C)$.
Let us set $x = \pi (z)$. Then, for any $W \in T (\partial C(M))$
\[ F_\theta (S , W)_z  = (\pi^* \theta )(W)_z \sigma (S)_z =
\frac{1}{n+2} \, \theta_x ((d_z \pi ) W_z ) = \]
\[ = g_{\theta , x} (T^\bot_{x} , (d_z \pi )W_z ) = 0 \]
as $(d_z \pi ) W_z$ is tangent to $\partial M$. Hence $S_z \in
{\rm Null}(j^* F_\theta )_z$ (and we may apply (\ref{e:incl})).
\par
Since $F_\theta (S,S) = 0$ and $S$ is tangent to $\partial C(M)$,
$F_\theta$ is indefinite on $T(\partial C(M))$. However (by the
first part of Theorem \ref{p:5}) $F_\theta$ is
non\-de\-ge\-ne\-ra\-te on $T(\pi^{-1}(\Omega ))$ hence $(j^*
F_\theta )_z$ has signature $(2n, 1)$ at each $z \in
\pi^{-1}(\Omega )$. $\square$
\begin{proposition} Let $M$ be a strictly pseudoconvex CR
manifold-with-boundary and $\theta$ a contact form with $G_\theta$
positive definite. Let $T$ be the characteristic direction of $d
\theta$. The property that $T \in T(\partial M)$ is not CR
invariant. If $T \in T(\partial M)$ and $\hat{T}$ is the
characteristic direction of $d \hat{\theta}$, where $\hat{\theta}
= e^{2u} \theta$ {\rm (}$u \in C^\infty (M)${\rm )}, then ${\rm
Sing}(\hat{T}^T ) = \emptyset$.
\end{proposition}
\noindent  {\em Proof}. Let us consider a local orthonormal (with
respect to $g_\theta$) frame of $T(\partial M)$ of the form $\{
E_1 , \cdots , E_{2n-1} , T \}$, so that $E_a \in H(M)$, $1 \leq a
\leq 2n-1$. Next, let us complete $\{ E_a \}$ to a local
orthonormal frame $\{ E_1 , \cdots , E_{2n} \}$ of $H(M)$ and set
$T_\alpha = (1/\sqrt{2}) (E_\alpha + i E_{\alpha + n})$, $1 \leq
\alpha \leq n$. Given another contact form $\hat{\theta} = e^{2u}
\theta$ ($u \in C^\infty (M)$) the characteristic direction of $d
\hat{\theta}$ is expressed by
\[ \hat{T} = e^{-2u} (T + i u^{\overline{\alpha}}
T_{\overline{\alpha}} - i u^\alpha T_\alpha ) = \]
\[ = e^{-2u} \{ T + \frac{i}{\sqrt{2}} (u^{\overline{\alpha}} -
u^\alpha ) E_\alpha + \frac{1}{\sqrt{2}} (u^{\overline{\alpha}} +
u^\alpha ) E_{\alpha + n } \} \] dove $u^\alpha = u_\alpha =
T_\alpha (u)$ (as $L_\theta (T_\alpha , T_{\overline{\beta}}) =
\delta_{\alpha\beta}$). Let $\xi$ be a unit normal on $\partial
M$. Then $\xi \in \{ \pm E_{2n} \}$ hence ${\rm Sing}(\hat{T}^T )
= {\rm Sing}(T) = \emptyset$. $\square$ \vskip 0.1in If $z \in
C(M)$ we denote by $\beta_z : T_{\pi (z)}(M) \to {\rm
Ker}(\sigma_z )$ the inverse of the ${\mathbb R}$-linear
isomorphism $d_z \pi : {\rm Ker}(\sigma_z ) \to T_{\pi (z)} (M)$.
It is an elementary matter that
\begin{lemma}
Given $v \in T_x (\partial M)$ its horizontal lift $\beta_z v$, $z
\in \pi^{-1} (x)$, is tangent to $\partial C(M)$.
\end{lemma}
Indeed, let $a : (-\epsilon , \epsilon ) \to
\partial M$ be a smooth curve such that $a(0) = x$ and $\dot{a}(0) = v$.
Let $X \in T(\partial M)$ be a tangent vector field such that $X_x
= v$. Let $a^\uparrow : (-\epsilon , \epsilon ) \to C(M)$ be the
unique horizontal lift of $a$, issuing at $z$. As $\pi (a^\uparrow
(t)) = a(t)$ one has $a^\uparrow (t) \in \partial C(M)$, $|t| <
\epsilon$. On the other hand $\dot{a}^\uparrow (0) \in {\rm
Ker}(\sigma_z )$ and it projects on $v$ hence
\[ T_z (\partial C(M)) \ni \dot{a}^\uparrow (0) = X^\uparrow_z =
\beta_z v. \] $\square$ \vskip 0.1in We set $T(\partial
M)^\uparrow = \{ \beta X : X \in T(\partial M) \}$ and ${\mathcal
V}_z = {\rm Ker}(d_z \pi )$, for $z \in \partial C(M)$. As
observed above, $\partial C(M)$ is tangent to the $S^1$-action
hence $\mathcal V$ is a smooth distribution on $\partial C(M)$.
\begin{lemma}
Let $M$ be a strictly pseudoconvex CR manifold-with-boundary. One
has the decomposition
\begin{equation}
T(\partial C(M)) = T(\partial M)^\uparrow \oplus {\mathcal V}
\label{e:old6}
\end{equation} Moreover, if $\partial M$ is tangent to the characteristic direction $T$ of
$d \theta$ then \begin{equation} T(\partial C(M))^\bot \subseteq
{\rm Ker}(\sigma ), \;\;\; (d \pi ) T(\partial C(M))^\bot
\subseteq H(M), \label{e:old7}
\end{equation}
\begin{equation}
{\rm Ker}(\sigma ) = T(\partial M)^\uparrow \oplus T(\partial
C(M))^\bot . \label{e:old8}
\end{equation} Here $T(\partial C(M))^\bot \to \partial C(M)$ is the
normal bundle of $j : \partial C(M) \hookrightarrow (C(M),
F_\theta )$. \label{p:6}
\end{lemma}
{\em Proof of Lemma} \ref{p:6}. Note that
\[ T(\partial M)^\uparrow \cap {\mathcal V} \subseteq {\rm Ker}(\sigma )
\cap {\rm Ker}(d \pi ) = (0), \] hence the sum $T(\partial
M)^\uparrow + {\mathcal V}$ is direct. The arguments preceding
Lemma \ref{p:6} show that $T(\partial M)^\uparrow \oplus {\mathcal
V} \subseteq T(\partial C(M))$. Viceversa, let $V \in T(\partial
C(M)) \subset T(C(M))$. Then (by the decomposition (\ref{e:old1}))
\begin{equation}
V = X^\uparrow + f \, S, \label{e:old9}
\end{equation}
for some $X \in T(M)$ and $f \in C^\infty (C(M))$. Then
\[ X_{\pi (z)} = (d_z \pi ) V_z \in T_{\pi (z)}(\partial M), \;\;\; z \in
\partial C(M), \] i.e. $X \in T(\partial M)$ and then $T(\partial C(M)) \subseteq
T(\partial M)^\uparrow \oplus {\mathcal V}$. To check
(\ref{e:old7}) let $V \in T(\partial C(M))^\bot \subset T(C(M))$
and use (\ref{e:old1}) to decompose as in (\ref{e:old9}). By
assumption $T \in T(\partial M)$ hence $T^\uparrow \in T(\partial
C(M))$ and then
\[ 0 = F_\theta (V , T^\uparrow ) = \tilde{G}((d \pi )V , (d \pi
)T^\uparrow ) + \theta ((d \pi ) T^\uparrow ) \sigma (V) = \]
\[ = \tilde{G}_\theta (X , T) + \frac{f}{n+2} = \frac{f}{n+2} \]
i.e. $f = 0$, or $V = X^\uparrow \in {\rm Ker}(\sigma )$. To check
the second statement in (\ref{e:old7}) let
\[ V \in T(\partial C(M))^\bot \subseteq {\rm Ker}(\sigma ) =
T(M)^\uparrow = H(M)^\uparrow \oplus ({\mathbb R} T)^\uparrow \]
i.e. $V = Y^\uparrow + f T^\uparrow$, for some $Y \in H(M)$.
Moreover $S \in {\rm Ker}(d \pi ) \subset T(\partial C(M))$, hence
$S$ and $V$ are orthogonal
\[ 0 = F_\theta (S , V) = \theta ((d \pi ) V) \sigma (S) =
\frac{f}{n+2} \] i.e. $f = 0$, or $V \in H(M)^\uparrow$.
(\ref{e:old7}) is proved and may be equivalently written
\[ T(\partial C(M))^\bot \subseteq H(M)^\uparrow . \]
When $T^\bot = 0$ the space $T(\partial C(M))$ is nondegenerate in
$(T(C(M)) , F_\theta )$ hence so does the perp space $T(\partial
C(M))^\bot$. Also
\[ T(C(M)) = T(\partial C(M)) \oplus T(\partial C(M))^\bot . \]
Let us prove (\ref{e:old8}). First
\[ T(\partial M)^\uparrow \cap T(\partial C(M))^\bot \subseteq T(\partial C(M))
\cap T(\partial C(M))^\bot = (0) \] hence the sum $T(\partial
M)^\uparrow + T(\partial C(M))^\bot$ is direct and (by
(\ref{e:old7}))
\begin{equation}
T(\partial M)^\uparrow \oplus T(\partial C(M))^\bot \subseteq {\rm
Ker}(\sigma ). \label{e:old10}
\end{equation} Finally (by (\ref{e:old6})) \[ {\rm
Ker}(\sigma ) \oplus {\rm Ker}(d \pi ) = T(C(M)) = T(\partial
C(M)) \oplus T(\partial C(M))^\bot =
\] \[ = T(\partial M)^\uparrow \oplus {\rm Ker}(d \pi ) \oplus
T(\partial C(M))^\bot \] and (\ref{e:old10}) yields
(\ref{e:old8}). $\square$ \vskip 0.1in
From now on we assume that
$\partial M$ is tangent to $T$. Then let us consider a local
orthonormal frame $\{ E_1 , \cdots , E_{2n-1} , T \}$ of
$T(\partial M)$, with respect to $i^* g_\theta$ (the first
fundamental form of $i :
\partial M \hookrightarrow M$), defined on some open set $U
\subseteq \partial M$. In particular $E_a \in H(M)$, $1 \leq a
\leq 2n-1$.
\begin{lemma}
Let $M$ be a strictly pseudoconvex CR
manifold-with-bo\-un\-da\-ry. Let $\theta$ be a contact form on
$M$ such that $G_\theta$ is positive definite and let $T$ be the
characteristic direction of $d \theta$. Assume that $\partial M$
is tangent to $T$. Then
\[ \{ E_1^\uparrow , \cdots , E_{2n-1}^\uparrow , T^\uparrow \pm
\frac{n+2}{2} \, S \} \] is a local orthonormal frame of
$T(\partial C(M))$, with respect to $j^* F_\theta$, defined on the
open set $\pi^{-1} (U) \subseteq \partial C(M)$. In particular
$T^\uparrow - ((n+2)/2)  S$ is a global timelike vector field on
$\partial C(M)$, i.e. $(\partial C(M) , j^* F_\theta )$ is a
spacetime. \label{p:7}
\end{lemma}
See also \cite{kn:BeEh}. The proof is straightforward.

\section{The geometry of the second fundamental form of the
boundaries}  As $(\partial C(M), j^* F_\theta )$ is a Lorentz
submanifold of $(C(M), F_\theta )$ we may write the Gauss equation
\[ \nabla^{C(M)}_X Y = \nabla^{\partial C(M)}_X Y + {\mathbb B}(X,Y), \]
for any $X, Y \in T(\partial C(M))$. Here $\nabla^{\partial C(M)}$
is the induced connection and ${\mathbb B}$ is the second
fundamental form of $j : \partial C(M) \hookrightarrow C(M)$. Cf.
e.g. \cite{kn:O'Ne}, p. 100. At this point, we wish to compute the
mean curvature vector of $j$
\[ {\mathbb H} = \frac{1}{2n+1} \; {\rm trace}_{j^* F_\theta}
({\mathbb B}). \] To this end it is convenient to use the local
frame in Proposition \ref{p:7}.
\begin{theorem} Let $M$ be a strictly pseudoconvex CR manifold-with-boundary, of
CR dimension $n$, and $\theta$ a contact form on $M$ such that
$G_\theta$ is positive definite. Assume that $\partial M$ is
tangent to the characteristic direction $T$ of $d \theta$. Let $\{
E_1 , \cdots , E_{2n-1} , T \}$ be a local $g_\theta$-orthonormal
frame of $T(\partial M)$ and $\xi$ a unit normal vector field on
$\partial M$, both defined on the open set $U \subseteq \partial
M$. Then the mean curvature vector ${\mathbb H}$ of the immersion
$j : \partial C(M) \hookrightarrow C(M)$ is given by
\begin{equation}\label{e:old11} {\mathbb H}_z = \frac{1}{2n+1} \,
\sum_{a=1}^{2n-1} g_\theta (\nabla_{E_a} E_a \, , \, \xi )_{\pi
(z)} \, \xi_{z}^\uparrow
\end{equation}
for any $z \in \pi^{-1}(U)$. Here $\nabla$ is the Tanaka-Webster
connection of $(M , \theta )$. In particular ${\mathbb H} =
(2n/(2n+1)) \, H^\uparrow$, where $H$ is the mean curvature vector
of the immersion $i : \partial M  \hookrightarrow M$. Therefore,
$\partial C(M)$ is minimal in $(C(M), F_\theta )$ if and only if
$\partial M$ is minimal in $(M , g_\theta )$. \label{t:1}
\end{theorem}
{\bf Example 5}. ${\mathbb R}^{2n}_+ \times {\mathbb R}$ is a
strictly pseudoconvex CR manifold (with the CR structure induced
from ${\mathbb H}_n$) whose boundary $N = \partial ({\mathbb
R}^{2n}_+ \times {\mathbb R})$ is tangent to $T =
\partial /\partial t$. The normal bundle of the boundary is the
span of $\xi =
\partial /\partial y^n - 2 x_n T$. By the Gauss formula, the second fundamental form
of the boundary is given by
\[ B(\frac{\partial}{\partial x^n} ,  \frac{\partial}{\partial x^n} ) = - 4
y_n \xi , \;\; B(\frac{\partial}{\partial x^\alpha} ,
\frac{\partial}{\partial x^n} ) = - 2 y_\alpha \xi , \;\;
B(\frac{\partial}{\partial x^\alpha} , \frac{\partial}{\partial
x^\beta} ) = 0, \]
\[ B(\frac{\partial}{\partial x^\alpha} ,
\frac{\partial}{\partial y^\beta}) = 0, \;\;
B(\frac{\partial}{\partial x^n} , \frac{\partial}{\partial
y^\beta} ) = 2 x_\beta \xi , \;\; B(\frac{\partial}{\partial
y^\alpha} , \frac{\partial}{\partial y^\beta}) = 0, \]
\[ B(\frac{\partial}{\partial x^\alpha} , T) = 0, \;\;
B(\frac{\partial}{\partial x^n} , T) = \xi , \;\;
B(\frac{\partial}{\partial y^\alpha} , T ) = 0, \;\; B(T , T) = 0.
\]
Here $1 \leq \alpha , \beta \leq 2n-1$. On the other hand, the
induced metric on $N$ is given by
\[ g : \left( \begin{array}{ccc} 2(\delta_{ij} + 2 y_i y_j ) & - 4
y_i x_\beta & - 2 y_i \\ - 4 x_\alpha y_j & 2(\delta_{\alpha\beta}
+ 2 x_\alpha x_\beta ) & 2 x_\alpha \\ - 2 y_j & 2 x_\beta & 1
\end{array} \right) \]
hence (by an argument similar to the proof of Lemma \ref{l:ex})
the corresponding cometric on $T^* (N)$ is given by
\begin{equation} \label{e:com} g^{-1} : \left( \begin{array}{ccc} \frac{1}{2}
\delta^{ij} & 0 & y^i \\ 0 & \frac{1}{2} \delta^{\alpha\beta} & -
x^\alpha \\ y^j & - x^\beta & 1 + 2 |x^\prime |^2 + 2 |y|^2
\end{array} \right)
\end{equation}
where $x^\prime = (x_1 , \cdots , x_{2n-1})$, $|x^\prime |^2 =
x_\alpha x^\alpha$ and $|y|^2 = y_j y^j$. Finally a calculation
(based on (\ref{e:com})) shows that $2n H = g^{ab} B(\partial_a ,
\partial_b ) = 0$, i.e. $N$ is minimal in $({\mathbb R}^{2n}_+
\times {\mathbb R} , g_{\theta_0})$. In particular (by Theorem
\ref{t:1}) $\partial C({\mathbb R}^{2n}_+ \times {\mathbb R})$ is
minimal in $(C({\mathbb R}^{2n}_+ \times {\mathbb R}) ,
F_{\theta_0})$. $\square$ \vskip 0.1in Let $\{ X_A : 1 \leq A \leq
2n+1 \}$ be a local $F_\theta$-orthonormal frame of $T(\partial
C(M))$, i.e. $F_\theta (X_A , X_B ) = \epsilon_A \delta_{AB}$,
with $\epsilon_1 = \cdots = \epsilon_{2n} = 1 = -\epsilon_{2n+1}$.
Then ${\mathbb H}$ is locally given by
\[ {\mathbb H} = \frac{1}{2n+1} \, \sum_A \epsilon_A {\mathbb B}(X_A ,
X_A ). \] \par {\em Proof of Theorem} \ref{t:1}. Using the local
frame furnished by Lemma \ref{p:7} we obtain
\begin{equation}
(2n+1) {\mathbb H} = \sum_{a=1}^{2n-1} {\mathbb B}(E_a^\uparrow ,
E_a^\uparrow ) + 2 (n+2) {\mathbb B}(T^\uparrow , S).
\label{e:old12}
\end{equation}
As a consequence of Lemma \ref{l:1} we have
\begin{equation}
\nabla^{C(M)}_{E_a^\uparrow} E_a^\uparrow = (\nabla_{E_a} E_a
)^\uparrow - \frac{n+2}{2} \; A(E_a , E_a ) S, \label{e:old13}
\end{equation}
\begin{equation}
\nabla^{C(M)}_{T^\uparrow} S = 0. \label{e:old14}
\end{equation}
The equation (\ref{e:old14}) implies ${\mathbb B}(T^\uparrow , S)
= 0$ (with the corresponding simplification of (\ref{e:old12})).
As $T \in T(\partial M)$ we have
\[ T(\partial M)^\bot \subseteq H(M). \]
We need the following
\begin{lemma} $\,$ \par\noindent
Assume that $\partial M$ is tangent to $T$. Let $T(\partial
M)^\bot \to \partial M$ be the normal bundle of the immersion $i :
\partial M \hookrightarrow M$. Then
\begin{equation}
[T(\partial M)^\bot ]^\uparrow = T(\partial C(M))^\bot .
\label{e:old15}
\end{equation}
\label{l:5}
\end{lemma}
{\em Proof of Lemma} \ref{l:5}. Let $\xi \in T(\partial M)^\bot$
and $V \in T(\partial C(M)) = T(\partial M)^\uparrow \oplus {\rm
Ker}(d \pi )$, i.e. $V = X^\uparrow + f \, S$. Let us set $X_H :=
X - \theta (X) T \in H(M)$. Then
\[ F_\theta (V , \xi^\uparrow ) = \tilde{G}_\theta (X , \xi ) + f
\; F_\theta (S, \xi^\uparrow ) = \]
\[ = G_\theta (X_H , \xi ) + f\; \theta (\xi ) \sigma (S) =
g_\theta (X_H , \xi ) = 0 \] because $X,T \in T(\partial M)$
implies $X_H \in T(\partial M)$. It follows that  \[ [T(\partial
M)^\bot ]^\uparrow \subseteq T(\partial C(M))^\bot . \] The
desired equality follows by inspecting dimensions. $\square$
\par
Let $\xi$ be a unit normal vector field on $\partial M$, defined
on the open set $U \subseteq N$. Then (by Lemma \ref{l:5})
$\xi^\uparrow$ is a unit normal vector field on $\partial C(M)$.
Then (by the Gauss equation and by (\ref{e:old13}))
\[ F_\theta ({\mathbb B}(E_a^\uparrow , E_a^\uparrow )\, , \,
\xi^\uparrow ) = F_\theta (\nabla^{C(M)}_{E_a^\uparrow}
E_a^\uparrow \, , \, \xi^\uparrow ) = F_\theta ((\nabla_{E_a} E_a
)^\uparrow \, , \, \xi^\uparrow ) = \]
\[ = \tilde{G}_\theta (\nabla_{E_a} E_a \, , \, \xi ) = g_\theta
(\nabla_{E_a} E_a \, , \, \xi ) \] which yields (\ref{e:old11}).
$\square$ \vskip 0.1in
The Levi-Civita connection
$\nabla^{g_\theta}$ of $(M , g_\theta )$ is related to the
Tanaka-Webster connection $\nabla$ of $(M , \theta )$ by
\begin{equation}
\nabla^{g_\theta}_X Y = \nabla_X Y + (\Omega (X,Y) - A(X,Y)) T +
\label{e:old16}
\end{equation}
\[ + \tau (X) \theta (Y) + \theta (X) J Y + \theta (Y) J X , \]
for any $X,Y \in T(M)$. Here $\Omega = - d \theta$. Cf. e.g.
\cite{kn:BaDr}, p. 238. Thus, for any $X, Y \in H(M)$
\[   \nabla^{g_\theta}_X Y = \nabla_X Y + (\Omega (X,Y) - A(X,Y)) T \]
and then
\[ \nabla^{g_\theta}_{E_a} E_a = \nabla_{E_a} E_a - A(E_a , E_a )
T \] implies  (as $g_\theta (T , \xi ) = 0$)
\[ (2n+1) {\mathbb H} = \sum_{a} g_\theta (\nabla^{g_\theta}_{E_a}
E_a \, , \, \xi ) \xi^\uparrow =    \]
\[ = \sum_{a} g_\theta (B(E_a , E_a ) \, , \, \xi ) \xi^\uparrow
= 2n g_\theta (H \, , \, \xi ) \xi^\uparrow \] because
$\nabla^{g_\theta}_T T = 0$ implies $B(T,T) = 0$. Here $B$ is the
second fundamental form of $i : \partial M \hookrightarrow M$ and
$H = (1/(2n)) \, {\rm trace}_{g_\theta} (B)$ is its mean curvature
vector. Then ${\mathbb H} = (2n/(2n+1)) \, H^\uparrow$. $\square$
\begin{theorem} Let $M$ be a strictly pseudoconvex CR
manifold-with-boundary and $\theta$ such that $T \in T(\partial
M)$. Then $\partial C(M)$ has nonumbilic points in $(C(M),
F_\theta )$. Moreover $\partial M$ is totally umbilical in $(M,
g_\theta )$ if and only if
\[ {\mathbb B}(X^\uparrow , Y^\uparrow ) = \frac{2n+1}{2n} \,
F_\theta (X^\uparrow , Y^\uparrow ) {\mathbb H}, \]
\[ {\mathbb B}(X^\uparrow , T^\uparrow ) = \{ (d \sigma
)(X^\uparrow , \xi^\uparrow ) + g_\theta (X , J \xi ) \}
\xi^\uparrow , \] for any $X, Y \in T(\partial M) \cap H(M)$.
\label{t:3}
\end{theorem}
\noindent {\em Proof}. By (\ref{e:old16}) and the Gauss formula
for the immersion $\partial M \hookrightarrow (M, g_\theta )$
\[ B(X,Y) = g_\theta (\nabla_X Y , \xi ) \xi , \;\;\;
B(X,T) = g_\theta (\tau X + J X , \xi ) \xi , \] for any $X,Y \in
T(\partial M) \cap H(M)$. Next, by Lemma \ref{l:1} and the Gauss
formula for the immersion $\partial C(M) \hookrightarrow (C(M),
F_\theta )$ \begin{equation}\label{e:sff1} {\mathbb B}(X^\uparrow
, Y^\uparrow ) = B(X,Y)^\uparrow , \end{equation}
\begin{equation}\label{e:sff2} {\mathbb B}(X^\uparrow , T^\uparrow ) =
B(X,T)^\uparrow + \{ (d \sigma )(X^\uparrow , \xi^\uparrow ) +
g_\theta (X , J \xi ) \} \xi^\uparrow , \end{equation}
\begin{equation}\label{e:sff3} {\mathbb B}(X^\uparrow , \hat{S}) = - g_\theta (X , J
\xi ) \xi^\uparrow , \;\;\; {\mathbb B}(T^\uparrow , \hat{S}) = 0.
\end{equation}
Note that $J \xi$ is tangent to $\partial M$. Assume that
${\mathbb B} = F_\theta \otimes {\mathbb H}$. Then (by
(\ref{e:sff3})) $J \xi$ is orthogonal to $\partial M$, hence $\xi
= 0$, a contradiction. The last statement in Theorem \ref{t:3}
follows from $B = g_\theta \otimes H$ and
(\ref{e:sff1})-(\ref{e:sff2}). $\square$

\section{Minimal submanifolds}
The purpose of this section to investigate minimal submanifolds in
the Heisenberg group ${\mathbb H}_n$. First, we establish the
relationship between the notion of $X$-minimality of N. Arcozzi \&
F. Ferrari, cf. (3) in \cite{kn:ArFe}, I. Birindelli \& E.
Lanconelli, cf. (3.23) in \cite{kn:BiLa}, and N. Garofalo \& S.D.
Pauls, cf. (2.5) in \cite{kn:GaPa} (see also \cite{kn:Pau}) and
minimality of an isometric immersion (between Riemannian
manifolds). Second, we prove the following
\begin{theorem}
Let $\Psi : N \to {\mathbb H}_n$ be an isometric immersion of a
$m$-dimensional Riemannian manifold $(N , g)$ into $({\mathbb H}_n
, g_{\theta_0})$. Then $\Psi$ is minimal if and only if
\begin{equation} \Delta \Psi = 2 J T^\bot
\label{e:delta}
\end{equation}
where $\Delta$ is the Laplace-Beltrami operator of $(N , g)$. In
particular, there are no minimal isometric immersions $\Psi$ of a
compact Riemannian manifold $N$ into the Heisenberg group such
that $T$ is tangent to $\Psi (N)$.
\label{t:min}
\end{theorem}
\noindent Compare to Theorem 6.2 and Corollaries 6.1 and 6.2 in
\cite{kn:DW}, p. 45-48. Let $M = {\mathbb H}_1$ be the lowest
dimensional Heisenberg group and $\varphi : {\mathbb H}_1 \to
{\mathbb R}$ a $C^2$ function. Let us set
\[ N = \{ x \in {\mathbb H}_1 : \varphi (x) = 0 \} \]
and assume there is an open neighborhood $O \supset N$ such that
\begin{equation}
|\nabla \varphi (x)| \geq \alpha
> 0, \;\;\; x \in O.
\label{e:old17}
\end{equation}
Here $\nabla \varphi$ is the Euclidean gradient of $\varphi$. Let
$(z, t)$ be the natural coordinates on ${\mathbb H}_1 = {\mathbb
C} \times {\mathbb R}$ and set $Z = Z_1 = \partial /\partial z + i
\, \overline{z} \, \partial /\partial t$ (the generator of
$T_{1,0}({\mathbb H}_1 )$). Let $\theta_0 = d t + i (z \, d
\overline{z} - \overline{z} \, d z)$ be the canonical contact form
on ${\mathbb H}_1$. Note that $L_{\theta_0} (Z , \overline{Z}) =
1$. The Tanaka-Webster connection of $({\mathbb H}_1 , \theta_0 )$
is given by
\[ \Gamma^A_{BC} = 0, \;\;\; A,B,C \in \{ 1 , \overline{1} , 0 \}
. \] Let us set $X_1 = \frac{1}{\sqrt{2}} \, (Z + \overline{Z})$
and  $X_2 = Y_1 = \frac{i}{\sqrt{2}} \, (Z - \overline{Z})$. We
shall prove the following
\begin{theorem} Let $N = \{ x \in {\mathbb H}_1 : \varphi (x) = 0
\}$ be a surface in ${\mathbb H}_1$ such that $(\ref{e:old17})$
holds. Assume that $N$ is tangent to the characteristic direction
$T = \partial /\partial t$ of $({\mathbb H}_1 , \theta_0 )$. Let
$\xi$ be a unit normal vector field on $N$. Then the mean
curvature vector of $N$ in $({\mathbb H}_1 , g_{\theta_0} )$ is
given by
\begin{equation}
H = - \frac{1}{2} \, \sum_{j=1}^2 X_j \left( \frac{X_j \varphi}{|X
\varphi |} \right) \; \xi . \label{e:old18}
\end{equation}
Here $| X \varphi |^2 = (X_1 \varphi )^2 + (X_2 \varphi )^2$ is
the $X$-gradient of $\varphi$.
\label{c:1}
\end{theorem}
{\em Proof of Theorem} \ref{c:1}. $T(N)$ is the span of $\{ E \, ,
\, T \}$ while $T(N)^\bot$ is the span of $\xi$, where
\[ E = \frac{1}{|X \varphi |} \{ (X_2 \varphi ) X_1 - (X_1
\varphi ) X_2 \} ,  \;\; \xi = \frac{1}{|X \varphi |} \{ (X_1
\varphi ) X_1 + (X_2 \varphi ) X_2 \} , \] so that $g_{\theta_0}
(E, E) = 1$ and $g_{\theta_0} (\xi , \xi ) = 1$. A calculation
(based on $\nabla_{X_j} X_k = 0$) leads to
\[ \nabla_{X_1} E = X_1 \left( \frac{X_2 \varphi}{|X \varphi |}
\right) X_1 - X_1 \left( \frac{X_1 \varphi}{|X \varphi |} \right)
X_2 \, , \]
\[ \nabla_{X_2} E = X_2 \left( \frac{X_2 \varphi}{|X \varphi |}
\right) X_1 - X_2 \left( \frac{X_1 \varphi}{|X \varphi |} \right)
X_2 \, , \] hence
\begin{equation}
\nabla_E E = \frac{1}{|X \varphi |} \left\{ \left[ (X_2 \varphi )
X_1 \left( \frac{X_2 \varphi}{|X \varphi |} \right) - (X_1 \varphi
) X_2 \left( \frac{X_2 \varphi}{| X \varphi |} \right) \right] X_1
+ \right. \label{e:old19}
\end{equation}
\[ + \left. \left[ (X_1 \varphi ) X_2 \left( \frac{X_1 \varphi}{|X
\varphi |} \right) - (X_2 \varphi ) X_1 \left( \frac{X_1
\varphi}{|X \varphi |} \right) \right] X_2 \right\} . \] Then (by
(\ref{e:old19}))
\begin{equation} g_{\theta_0} (\nabla_E E , \xi ) = - \sum_{j=1}^2 X_j
\left( \frac{X_j \varphi}{|X \varphi |} \right) + \label{e:old20}
\end{equation}
\[ + \frac{1}{|X \varphi |^2} \left\{ (X_1 \varphi )^2 \; X_2 \left(
\frac{X_2 \varphi}{|X \varphi |} \right) + (X_1 \varphi )^2 \; X_1
\left( \frac{X_1 \varphi}{| X \varphi |} \right) + \right. \]
\[ + \left. (X_1 \varphi )(X_2 \varphi ) \; X_1 \left( \frac{X_2
\varphi}{|X \varphi |} \right) + (X_1 \varphi )(X_2 \varphi ) \;
X_2 \left( \frac{X_1 \varphi}{|X \varphi |} \right) \right\} . \]
Using the identity
\[ |X \varphi | X_j (|X \varphi |) = (X_1 \varphi ) \; X_j X_1
\varphi + (X_2 \varphi ) \; X_j X_2 \varphi \] one may show that
the second term in the right hand member of (\ref{e:old20}) is $|X
\varphi |^{-4}$ times
\[ (X_2 \varphi )^2 \{ (X_2 X_2 \varphi ) |X \varphi | - (X_2
\varphi ) X_2 (|X \varphi |) \} + \]
\[ + (X_1 \varphi )^2 \{ (X_1 X_1 \varphi )|X \varphi | -
(X_1 \varphi ) X_1 (|X \varphi |) \} + \]
\[ + (X_1 \varphi )(X_2 \varphi ) \{ (X_1 X_2 \varphi ) |X \varphi
| - (X_2 \varphi ) X_1 (|X \varphi |) \} + \]
\[ + (X_1 \varphi )(X_2 \varphi ) \{ (X_2 X_1 \varphi ) |X \varphi
| - (X_1 \varphi ) X_2 (|X \varphi |) \} = \]
\[ = - \{ (X_1 \varphi ) X_1 (|X \varphi |) + (X_2 \varphi ) X_2 (|X \varphi |) \}
\{ (X_1 \varphi )^2 + (X_2 \varphi )^2 \} + \]
\[ + |X \varphi | \{ (X_1 \varphi )^2 X_1 X_1 \varphi +
2 (X_1 \varphi )(X_2 \varphi ) X_1 X_2 \varphi +(X_2 \varphi )^2
X_2 X_2 \varphi \} \] (as $[X_1 , X_2 ] = - 2 \; T$ and $T(\varphi
) = 0$) or
\[ - (X_2 \varphi ) |X \varphi | \{ (X_1 \varphi ) X_2 X_1 \varphi
+ (X_2 \varphi ) X_2 X_2 \varphi \} - \]
\[ - (X_1 \varphi ) |X \varphi | \{ (X_1 \varphi ) X_1
X_1 \varphi + (X_2 \varphi ) X_1 X_2 \varphi \} + \]
\[ + |X \varphi | \{ (X_1
\varphi )^2 X_1 X_1 \varphi + 2 (X_1 \varphi ) (X_2 \varphi ) X_1
X_2 \varphi + (X_2 \varphi )^2 X_2 X_2 \varphi \} = 0 \] hence
(\ref{e:old20}) leads to (\ref{e:old18}). $\square$ \vskip 0.1in
Let us prove Theorem \ref{t:min}. Let $(x^1 , \cdots , x^{2n} ,
x^0 )$ be the Cartesian coordinates on ${\mathbb R}^{2n+1}$ and
$(U, u^1 , \cdots , u^m )$ a local coordinate system on $N$. Let
$H(\Psi )$ be the mean curvature vector of $\Psi : N \to {\mathbb
H}_n$. Then $H(\Psi ) = H^A \partial_A$, where $\partial_A$ is
short for $\partial /\partial x^A$. Let $g_0 = g_{\theta_0}$ be
the Webster metric of $({\mathbb H}_n , \theta_0 )$ and $D^0$ the
Levi-Civita connection of $({\mathbb H}_n , g_0 )$. We set
$B_\alpha^A = \partial \Psi^A /\partial u^\alpha$, so that $\Psi_*
(\partial /\partial u^\alpha ) = B^A_\alpha \partial_A$. Let $\{
E_1 , \cdots , E_m \}$ be a local orthonormal (with respect to
$g$) frame of $T(N)$, defined on $U$. Then $E_\alpha =
E_\alpha^\beta \partial /\partial u^\beta$. Taking into account
that $E_\alpha^\beta B^A_\beta = E_\alpha (\Psi^A )$, the Gauss
formula of $\Psi$
\[ D^0_{E_\alpha} E_\beta = \Psi_* D_{E_\alpha} E_\beta + B(E_\alpha ,
E_\beta ) \] may be written
\[ \{ E_\alpha (E_\beta \Psi^A )  - (D_{E_\alpha} E_\beta
)(\Psi^A ) \} \partial_A   = \] \[ = B(E_\alpha , E_\beta ) -
E_\alpha (\Psi^A ) E_\beta (\Psi^B ) D^0_{\partial_A} \partial_B
\, . \] Here $D$ is the Levi-Civita connection of $(N , g)$ and
$B$ is the second fundamental form of $\Psi$. Contraction of
$\alpha$ and $\beta$ gives
\begin{equation}
(\Delta \Psi^A )\partial_A = m H(\Psi ) - \sum_{\alpha = 1}^m
E_\alpha (\Psi^A ) E_\alpha (\Psi^B ) D^0_{\partial_A} \partial_B
\label{e:min33}
\end{equation}
Since
\begin{equation}
\partial_j = \frac{\partial}{\partial x^j}
= Z_j + \overline{Z}_j - 2 y^j T, \;\; \partial_{j+n} =
\frac{\partial}{\partial y^j} = i(Z_j - \overline{Z}_j ) + 2 x^j
T,
\label{e:partials}
\end{equation} it follows that the
Tanaka-Webster connection of $({\mathbb H}_n , \theta_0 )$
satisfies
\[ \nabla_{\partial_j} \partial_k = \nabla_{\partial_{j+n}}
\partial_{k+n} = 0, \]
\begin{equation} \nabla_{\partial_j} \partial_{k+n} = -
\nabla_{\partial_{j+n}}
\partial_k = 2 \delta_{jk} T,
\label{e:min34}
\end{equation}
\[ \nabla_{\partial_A} T = \nabla_T \partial_B = 0. \]
Let $J$ be the complex structure in $H({\mathbb H}_n )$, extended
to a $(1,1)$-tensor field on ${\mathbb H}_n$ by requesting that $J
T = 0$. Using $D^0 = \nabla - (d \theta_0 ) \otimes T + 2
(\theta_0 \odot J)$ it follows that
\[ E_\alpha (\Psi^A ) E_\beta (\Psi^B ) D^0_{\partial_A}
\partial_B = E_\alpha (\Psi^A ) E_\beta (\Psi^B )
\nabla_{\partial_A} \partial_B - \] \[ - (d \theta )(E_\alpha ,
E_\beta )T + \theta (E_\alpha ) J \Psi_* E_\beta + \theta (E_\beta
) J \Psi_* E_\alpha \] where $\theta = \Psi^* \theta_0$. On the
other hand, by (\ref{e:min34})
\[ E_\alpha (\Psi^A ) E_\beta (\Psi^B )
\nabla_{\partial_A} \partial_B = \] \[ = 2 \sum_{j=1}^n \{
E_\alpha (\Psi^j ) E_\beta (\Psi^{j+n}) - E_\alpha (\Psi^{j+n})
E_\beta (\Psi^{j}) \} T . \] Also $\sum_\alpha \theta (E_\alpha )
= \sum_\alpha g_0 (T , \Psi_* E_\alpha )  = \sum_\alpha g(T^T ,
E_\alpha )$ hence
\[ \sum_\alpha \theta (E_\alpha ) J \Psi_* E_\alpha = J \Psi_* T^T
= - J T^\bot , \] so that (\ref{e:min33}) becomes $m H(\Psi ) =
\Delta \Psi - 2 J T^\bot$ (yielding (\ref{e:delta})). $\square$
\vskip 0.1in Our Theorem \ref{c:1} demonstrates that the Webster
metric is the "correct" choice of ambient metric. Nevertheless,
even the geometry of a hyperplane in $({\mathbb H}_n \, , \, g_0
)$ turns out to be rather involved. In the sequel, we work out
explicitly the case of $\{ z \in {\mathbb H}_n : t = 0 \}$. \vskip
0.1in {\bf Example 1.} {\rm ({\em continued}) Let $\Psi :
\partial {\mathbb H}_n^+ \to {\mathbb H}_n^+$ be the inclusion and
$g = \Psi^* g_0$ (the first fundamental form of $\Psi$). Let
$\Delta$ be the Laplace-Beltrami operator of $(\partial {\mathbb
H}^n_+ , g)$. We may state
\begin{proposition} The coordinate functions $z^j$  on
$\partial {\mathbb H}_n^+ \approx {\mathbb C}^n$ satisfy $\Delta
z^j = 2\, z^j /(1 + 2 |z|^2 )$. Consequently the boundary of
$({\mathbb H}^+_n , g_0 )$ is mi\-nimal. \label{p:ex1}
\end{proposition} \noindent Note that
\[ \theta_0 (\partial_i ) = - 2 y_i \, , \;\;\; \theta_0 (\partial_{i+n}) = 2 x_i \, , \]
\[ (d \theta_0 )(\partial_i , \partial_j ) = (d \theta_0
)(\partial_{i+n} , \partial_{j+n}) = 0, \;\; (d \theta_0
)(\partial_i , \partial_{j+n}) = 2 \delta_{ij} \, , \] \[ J
\partial_j = \partial_{j+n} - 2 x_j T , \;\;\; J \partial_{j+n} = - \partial_j
- 2 y_j T. \] Then by (\ref{e:old16}) (with $\tau = 0$) and by
(\ref{e:min34}) it follows that \[ D^0_{\partial_i}
\partial_j = - 2 (y_i \delta^k_j + y_j \delta_i^k ) \frac{\partial}{\partial y^k} +
4 (y_i x_j + y_j x_i ) T, \]
\[ D^0_{\partial_i} \partial_{j+n} = 2 (y_i
\frac{\partial}{\partial x^j} + x_j \frac{\partial}{\partial y^i}
) + 4(y_i y_j - x_i x_j ) T, \]
\[ D^0_{\partial_{i+n}} \partial_{j+n} = -2(x_i \delta_j^k + x_j
\delta_i^k ) \frac{\partial}{\partial x^k} - 4 (x_i y_j + x_j y_i
) T. \]
Next, we shall need the Gauss formula
\[ D^0_{\partial_a} \partial_b = D_{\partial_a} \partial_b +
B(\partial_a , \partial_b ) , \] where $D$ is the Levi-Civita
connection of $(\partial {\mathbb H}_n^+ , g)$. We obtain
\[ B(\partial_i , \partial_j ) = 4c(y_i x_j + y_j x_i ) \xi , \]
\begin{equation}
B(\partial_i , \partial_{j+n}) = 4c(y_i y_j - x_i x_j ) \xi ,
\label{e:B}
\end{equation}
\[ B(\partial_{i+n} , \partial_{j+n}) = - 4 c (x_i y_j + x_j y_i )
\xi , \] hence $\Psi$ is not totally geodesic, and if
$D_{\partial_a} \partial_b = \Gamma^c_{ab} \partial_c$ then
\[ \Gamma^k_{ij}  = - 4 c
(y_i x_j + y_j x_i ) y^k \, , \;\; \Gamma^{k + n}_{i + n \; j + n}
= - 4 c (x_i y_j + x_j y_i ) x^k \, , \]
\[ \Gamma^{k+n}_{ij} = 4 c
(y_i x_j + y_j x_i ) x^k - 2 (y_i \delta_j^k + y_j \delta_i^k ) ,
\]
\begin{equation}
\Gamma^k_{i \; j+n} = 2 y_i \delta_j^k - 4 c (y_i y_j - x_i x_j )
y^k \, ,
\label{e:G}
\end{equation}
\[ \Gamma^{k+n}_{i \; j+n} = 2 x_j \delta_i^k + 4 c (y_i y_j - x_i
x_j ) x^k \, , \]
\[ \Gamma^k_{i+n \; j+n} = - 2(x_i \delta_j^k + x_j \delta_i^k ) +
4 c (x_i y_j + x_j y_i ) y^k \, . \] We need the following
\begin{lemma} The local coefficients of the cometric $g^{-1}$ on
$T^* (\partial {\mathbb H}_n^+ )$ are given by
\begin{equation}
\label{e:gmeno} g^{-1} : \left(
\begin{array}{cc} \frac{1}{2} \delta^{ij} - c y^i y^j & c y^i x^j
\\ c x^i y^j & \frac{1}{2} \delta^{ij} - c x^i x^j \end{array}
\right) .
\end{equation} Consequently
\[ \Delta u = \frac{1}{2} \, \Delta_0 u + 2 c \frac{\partial
u}{\partial r} - c \{ y^i y^j \frac{\partial^2 u}{\partial x^i
\partial x^j} - 2 y^i x^j \frac{\partial^2 u}{\partial x^i
\partial y^j} + x^i x^j \frac{\partial^2 u}{\partial y^i \partial
y^j} \} \] for any $u \in C^2 (\partial {\mathbb H}_n^+ )$, where
$\Delta_0$ is the ordinary Laplacian on ${\mathbb R}^{2n}$ and
$\partial /\partial r$ is the radial vector field $x^j (\partial
/\partial x^j ) + y^j (\partial /\partial y^j )$. \label{l:ex}
\end{lemma}
\noindent By Lemma \ref{l:ex} it follows that $\Delta x^j = 2 c
x^j$ and $\Delta y^j = 2 c y^j$, hence the first statement in
Proposition \ref{p:ex1}. On the other hand $T^\bot = c \xi$
implies $J T^\bot = c \; \partial /\partial r$ hence (by Theorem
\ref{t:min}) $H(\Psi ) = 0$. Note that the mean curvature vector
may be also computed from $2n H(\Psi ) = g^{ab} B(\partial_a ,
\partial_b )$ by (\ref{e:B}) and (\ref{e:gmeno}).
\par
It remains that we prove Lemma \ref{l:ex}. The first statement is
elementary yet rather involved. The identities $g_{ac} g^{cb} =
\delta_a^b$ may be written \begin{equation} \begin{cases} 2
(\delta_{ij} + 2 y_i y_j ) g^{jk} - 4 y_i x_j g^{j+n,k} =
\delta^k_i \, , \cr  2(\delta_{ij} + 2 y_i y_j ) g^{j,k+n} - 4 y_i
x_j g^{j+n,k+n} = 0, \cr  - 4 x_i y_j g^{jk} + 2(\delta_{ij} + 2
x_i x_j ) g^{j+n,k} = 0, \cr - 4 x_i y_j g^{j,k+n} + 2
(\delta_{ij} + 2 x_i x_j ) g^{j+n,k+n}= \delta^k_i \, . \cr
\end{cases}
\label{e:expl}
\end{equation}
 Contraction of the first two equations (respectively
of the last two equations) by $y^i$ (respectively by $x^i$) gives
\[ (1 + 2 |y|^2 ) y_j g^{jk} - 2 |y|^2 x_j g^{j+n,k} = y^k \, ,
\]
\[ (1+ 2|y|^2 ) y_j g^{j,k+n} - 2 |y|^2 x_j g^{j+n,k+n} = 0, \]
\[ 2 |x|^2 y_j g^{jk} - (1 + 2 |x|^2 ) x_j g^{j+n,k} = 0, \]
\[ - 2 |x|^2 y_j g^{j,k+n} + (1 + 2 |x|^2 ) x_j g^{j+n,k+n} = x^k
\, , \] where from
\[ y_j g^{jk} = \frac{c}{2} (1 + 2 |x|^2 ) y^k \, , \;\; x_j
g^{j+n,k} = c |x|^2 y^k \, , \]
\[ y_j g^{j,k+n} = c |y|^2 x^k \, , \;\; x_j g^{j+n,k+n} = \frac{c}{2} (1 + 2 |y|^2 ) x^k
\, , \] and substitution back into (\ref{e:expl}) yields
(\ref{e:gmeno}). To compute the Laplacian
\[ \Delta u = \frac{\partial}{\partial x^a} \left( g^{ab}
\frac{\partial u}{\partial x^b}\right) + g^{ab}
\frac{\partial}{\partial x^a} \left( \log \sqrt{G} \right)
\frac{\partial u}{\partial x^b} \] (with $G = \det [g_{ab}]$) we
recall that $\partial (\log \sqrt{G} )/\partial x^a =
\Gamma^b_{ba}$ hence (by (\ref{e:G}))
\[ \frac{\partial}{\partial x^a} \left( \log \sqrt{G} \right) = 2
c x_a \, , \;\;\; 1 \leq a \leq 2n . \] Then (\ref{e:gmeno})
yields the result.} $\square$

\section{The CR Yamabe problem}
Let $M$ be a compact strictly pseudoconvex CR
manifold-with-bo\-un\-da\-ry, of CR dimension $n$, and $\theta$ a
contact form on $M$ with $G_\theta$ positive definite. Let us
assume that $\partial M$ is tangent to the characteristic
direction $T$ of $d \theta$.
\begin{lemma} Let us
set $p = 2 + 2/n$ and $f = (p-2) \log u$, with $u \in C^\infty
(M)$, $u > 0$. If $\hat{\theta} = e^f \theta$ then $\partial C(M)$
is minimal in $(C(M), F_{\hat{\theta}})$ if and only if
\begin{equation}
\label{e:new29}
\frac{\partial (u \circ \pi )}{\partial \eta} - n\;  F_\theta
({\mathbb H}, \eta) \; u \circ \pi = 0 \;\;\; {\rm on} \;\;
\partial C(M),
\end{equation} where $\eta$ and ${\mathbb
H}$ are respectively an outward unit normal and the mean curvature
vector of the immersion $\partial C(M) \hookrightarrow (C(M),
F_\theta )$. In particular, if $\xi$ and $H$ are an outward unit
normal and the mean curvature vector of the immersion $\partial M
\hookrightarrow (M , g_\theta )$ then {\rm (\ref{e:new29})}
projects to
\begin{equation}
\label{e:new30}
\frac{\partial u}{\partial \xi} -
\frac{2n^2}{2n+1} \; g_\theta (H , \xi ) \; u = 0 \;\;\; {\rm on}
\;\; \partial M.
\end{equation}
 \label{l:3}
\end{lemma}
\noindent The first statement in Lemma \ref{l:3} is of course well
known in conformal geometry. We give a brief proof for the
convenience of the reader. If $\hat{\theta} = e^f \theta$ the
corresponding Fefferman metric is $F_{\hat{\theta}} = e^{f \circ
\pi} F_\theta$ hence the Levi-Civita connections $\hat{D}$ and $D$
(of $F_{\hat{\theta}}$ and $F_\theta$, respectively) are related
by
\begin{equation}
\hat{D}_V W = D_V W + \frac{1}{2} \{ V(f) W + W(f) V - F_\theta
(V,W) D (f\circ \pi ) \} ,
\label{e:29}
\end{equation}
for any $V, W \in T(C(M))$, where $D (f\circ \pi )$ is the
gradient of $f\circ \pi$ with respect to $F_\theta$. Our
assumption $T \in T(\partial M)$ and Proposition \ref{p:5} imply
that $T(\partial C(M))$ is nondegenerate in $T(C(M))$ with respect
to $F_\theta$, hence with respect to $F_{\hat{\theta}}$ as well.
Let $\mathbb B$ and $\hat{\mathbb B}$ be the second fundamental
forms of the immersions $\partial C(M) \hookrightarrow (C(M),
F_\theta )$ and $\partial C(M) \hookrightarrow (C(M),
F_{\hat{\theta}})$. Then (by (\ref{e:29}) and the Gauss formula)
\begin{equation}
\hat{\mathbb B} = {\mathbb B} - \frac{1}{2} \; F_\theta \otimes
(D(f \circ \pi ))^\bot \, . \label{e:30}
\end{equation}
Taking traces in (\ref{e:30}) shows that the mean curvature
vectors of the two immersions are related by $\hat{\mathbb H} =
e^{-f} \{ {\mathbb H} - \frac{1}{2} (D(f \circ \pi ))^\bot \}$
hence $\partial C(M)$ is minimal in $(C(M) , F_{\hat{\theta}})$ if
and only if ${\mathbb H} = (1/2) (D(f \circ \pi ))^\bot$ and
(\ref{e:new29}) is proved. Let $\xi$ be an outward unit normal on
$\partial M$ in $(M , g_\theta )$. Then $\eta = \xi^\uparrow$ is
an outward unit normal on $\partial C(M)$ in $(C(M), F_{\theta})$.
Then (by Theorem \ref{t:1}) the mean curvatures of $\partial M
\hookrightarrow (M , g_\theta )$ and $\partial C(M)
\hookrightarrow (C(M), F_\theta )$ are related by
\[ F_\theta ({\mathbb H}, \eta ) = \frac{2n}{2n+1} \; g_\theta (H , \xi ) \circ
\pi \] hence (\ref{e:new29}) projects on $M$ to give
(\ref{e:new30}). $\square$ \vskip 0.1in We may consider the
problem
\begin{equation}
- b_n \, \Delta_b u + \rho \; u = \lambda \; u^{p-1} \;\;\; {\rm
in} \;\; M,
\label{e:33}
\end{equation}
\begin{equation}
\frac{\partial u}{\partial \xi} - \frac{2n^2}{2n+1} \; \mu_\theta
\; u = 0 \;\;\; {\rm on} \;\; \partial M,
\label{e:34}
\end{equation}
(the {\em CR Yamabe problem} on a CR manifold-with-boundary) where
\[ \Delta_b u = {\rm
div}(\nabla^H u), \;\;\; u \in C^2 (M), \] is the {\em
sublaplacian} of $(M , \theta )$, $b_n = 2 + 2/n$, $\lambda$ is a
constant, and $\mu_\theta = g_\theta (H , \xi ) \in \{ \pm \| H \|
\}$. Also $\nabla u$ is the gradient of $u$ with respect to
$g_\theta$ and $\nabla^H u = \pi_H \nabla u$ (the {\em horizontal
gradient}) where $\pi_H : T(M) \to H(M)$ is the projection
associated with the direct sum decomposition $T(M) = H(M) \oplus
{\mathbb R} T$. The divergence operator is meant with respect to
the volume form $\omega = \theta \wedge (d \theta )^n$. The
problem (\ref{e:33})-(\ref{e:34}) is a nonlinear subelliptic
problem of variational origin. Indeed, we may state
\begin{theorem}
Let us set
\[ A_\theta (u) = \int_M \{ b_n \| \nabla^H u \|^2 + \rho \; u^2
\} \omega - a_n \int_{\partial M} \mu_\theta \; u^2 \; d \sigma \,
, \]
\[ B_\theta (u) = \int_M |u|^p \omega \, , \]
where $\sigma = {\rm vol}(i^* g_\theta )$, the canonical volume
form associated with the induced metric $i^* g_\theta$ on
$\partial M$, and $a_n = 2^{n+2} \, (n + 1)! \, n/(2n+1)$.
Moreover, let
\[ Q_\theta (u) = \frac{A_\theta (u)}{B_\theta (u)}\, , \;\;\;
Q(M) = \inf \{ Q_\theta (u) : u \in C^\infty (M), \;\; u > 0 \} .
\] If $u \in C^\infty (M)$ is a positive function such that $Q_\theta (u) =
Q(M)$ then $u$ is a solution to {\rm (\ref{e:33})-(\ref{e:34})}
with $\lambda = (p/2) Q(M)$, a CR invariant of $M$.
\label{p:11}
\end{theorem}
\noindent {\em Proof}. If $\{ T_\alpha \}$ is a local frame of
$T_{1,0}(M)$ then the horizontal gradient is expressed by
$\nabla^H u = u^\alpha T_\alpha + u^{\overline{\alpha}}
T_{\overline{\alpha}}$, where $u^\alpha =
g^{\alpha\overline{\beta}} u_{\overline{\beta}}$ and
$u_{\overline{\beta}} = T_{\overline{\beta}}(u)$, hence $\|
\nabla^H u \|^2 = 2 u_\alpha u^\alpha$. Then
\[ \frac{d}{d t} \{ A_\theta (u + t h) \}_{t=0} = 2 \int_M \{ b_n
(u^\alpha h_\alpha + u_\alpha h^\alpha ) + \rho u h \} \omega - 2
a_n \int_{\partial M} \mu_\theta \, u \, h \; d \sigma , \] for
any $h \in C^2 ({\rm Int}(M)) \cap C^1 (M)$ (where ${\rm Int}(M) =
M \setminus \partial M$). On the other hand
\[ \int_M u^\alpha h_\alpha \, \omega = \int_M \{ T_\alpha
(u^\alpha h) - h T_\alpha (u^\alpha ) \} \omega = \] \[ = \int_M
{\rm div}(h u^\alpha T_\alpha ) \omega - \int_M \{ T_\alpha
(u^\alpha ) + u^\alpha {\rm div}(T_\alpha ) \} h \, \omega . \]
Note that ${\rm div}(T_\alpha ) = \Gamma^\beta_{\beta\alpha}$
hence $T_\alpha (u^\alpha ) + u^\alpha {\rm div}(T_\alpha ) =
{u^\alpha}_\alpha$, where ${u^\alpha}_\beta =
g^{\alpha\overline{\gamma}} u_{\overline{\gamma}\beta}$ and
$u_{\alpha\overline{\beta}} = (\nabla^2 u)(T_\alpha ,
T_{\overline{\beta}})$. The complex Hessian is meant with respect
to the Tanaka-Webster connection i.e. \[ (\nabla^2 u)(X,Y) =
(\nabla_X d u )Y = X(Y(u)) - (\nabla_X Y)(u), \] for any $X, Y \in
{\mathcal X}(M)$. Note that $\omega = c_n \; d \, {\rm
vol}(g_\theta )$ (with $c_n = 2^n n!$). Then (by Green's lemma)
\[ \int_M u^\alpha h_\alpha \, \omega = c_n \int_{\partial M} h\,
u^\alpha g_\theta (T_\alpha , \xi ) d \sigma - \int_M
{u^\alpha}_\alpha \, h \, \omega . \] As the sublaplacian is
locally given by \[ \Delta_b u = {u^\alpha}_\alpha +
{u^{\overline{\alpha}}}_{\overline{\alpha}}
\] we may conclude that
\begin{equation}
\frac{d}{d t} \{ A_\theta (u + t h) \}_{t=0} = 2 \int_{M} (- b_n
\, \Delta_b u + \rho \, u) h \, \omega +
\label{e:35}
\end{equation}
\[ + 2 \int_{\partial M} [b_n c_n g_\theta (\nabla^H u , \xi ) -
a_n \mu_\theta u ] h \, d \sigma . \] Also
\begin{equation}
\frac{d}{d t} \{B_\theta (u + t h)\}_{t=0} = p \int_M u^{1 + 2/n}
h \, \omega . \label{e:36}
\end{equation}
As $T \in T(\partial M)$ one has $\xi \in H(M)$ hence $g_\theta
(\nabla^H u , \xi ) = \xi (u)$ (also denoted by $\partial
u/\partial \xi$). If $u$ achieves $Q(M)$
\[ \frac{d}{d t} \{ Q_\theta (u + t h)\}_{t=0} = 0 \]
hence
\[ 2 \int_M (-b_n \, \Delta_b u + \rho \, u) h \, \omega + 2
\int_{\partial M} [b_n c_n \, \xi (u) - a_n \mu_\theta u ] h \, d
\sigma - \] \[ - p \; Q_\theta (u) \, \int_M u^{1 + 2/n} h \,
\omega = 0. \] In particular this holds for $\left.
h\right|_{\partial M} = 0$ hence  \[ -b_n \, \Delta_b u + \rho u =
(p/2) Q(M) u^{1 + 2/n}
\] and going back to arbitrary $h$
\[ \frac{\partial u}{\partial \xi} - \frac{a_n}{b_n c_n} \, \mu_\theta \, u = 0
\;\;\; {\rm on} \;\; \partial M \] which is (\ref{e:34}) because
$a_n /(b_n c_n ) = 2n^2 /(2n+1)$. The proof that $Q(M)$ is a CR
invariant is similar to the arguments in \cite{kn:JL2}, p.
174-175. Let $E^+ \to M$ be the ${\mathbb R}_+$-bundle spanned by
$\theta$ and let us set
\[ E^\alpha_x = \{ \nu : E^+_x \to {\mathbb R} : \nu (t \theta_x )
= t^{-\alpha} \nu (\theta_x ), \; {\rm for \; all} \; t > 0\} ,
\;\; (\alpha > 0) \] for any $x \in M$. Then $(\nu_\theta )_x (t
\theta_x ) = 1/t$ defines a global frame $\{ \nu_\theta \}$ of
$E^1 \to M$ (and of course $\{ \nu_\theta^\alpha \}$ is a global
frame of $E^\alpha \to M$). We need the {\em CR invariant
sublaplacian}
\[ L : \Gamma^\infty (E^{n/2}) \to \Gamma^\infty (E^{1 + n/2}),
\;\; L (u \, \nu_\theta^{n/2}) = (- b_n \, \Delta_b u + \rho \, u)
\nu_\theta^{1 + n/2} . \] By definition $\int_M u \,
\nu_\theta^{n+1} = \int_M u \, \omega$. A section $s = u
\nu^\alpha_\theta$ in $E^\alpha$ is {\em positive} if $u > 0$.
Finally, the fact that $Q(M)$ is a CR invariant follows from
\begin{equation} \label{e:37}
Q(M) = \inf \{ \int_M (L s) \otimes s \; : \; s \in \Gamma^\infty
(E^{n/2})
\end{equation} \[  {\rm  a \;
positive \; section \; such \; that} \;\; \int_M s^p = 1 \} .
\]
The identity (\ref{e:37}) follows from the fact that the sets $\{
A_\theta (u) : B_\theta (u) = 1 , u > 0 \}$ and $\{ A_\theta
(u)/B_\theta (u) : u > 0  \}$ coincide and from the calculation
\[ \int_M (L s) \otimes s = \int_M (-b_n \, u \Delta_b u +
\rho u^2 ) \omega ,  \]
\[ \int_M u (\Delta_b u) \; \omega = \int_M \{ {\rm div}(u \nabla^H u)
- \| \nabla^H u \|^2 \} \omega = \]
\[ = c_n \int_{\partial M} u \, \frac{\partial u}{\partial \xi} \,
d \sigma - \int_M \| \nabla^H u \|^2  \omega , \] hence (by
(\ref{e:34})) $\int_M (L s) \otimes s = A_\theta (u)$, for any $s
= u \nu_\theta^{n/2} \in \Gamma^\infty (E^{n/2})$.

\section{Minimal surfaces in ${\mathbb H}_n$}
Let $(N, g)$ be a $2$-dimensional Riemannian manifold and $\Psi :
N \to {\mathbb H}_n$ a minimal isometric immersion of $(N , g)$
into $({\mathbb H}_n , g_0 )$. Let $(U , z = x + i y)$ be
isothermal local coordinates on $N$, i.e. locally \[ g = 2 E (d
x^2 + d y^2 ), \] for some $E \in C^\infty (U)$, $E > 0$. As well
known the Laplace-Beltrami operator of $(N , g)$ is locally given
by \[ \Delta u = \frac{2}{E} \; \frac{\partial^2 u}{\partial z
\partial \overline{z}} \, , \;\;\; u \in C^2 (N). \]
Let us set
$F^j = \Psi^j + i \Psi^{j+n}$, $1 \leq j \leq n$, and $f =
\Psi^0$. Also, we consider $K : U \to {\mathbb C}$ given by
\[ K = \frac{\partial f}{\partial z} + i \sum_{j=1}^n (F^j
\frac{\partial \overline{F}^j}{\partial z} - \overline{F}^j
\frac{\partial F^j}{\partial z}) . \]
\begin{lemma} The normal component of the characteristic vector
field $T = \partial /\partial t$ of $d \theta_0$ is locally given
by
\begin{equation} T^\bot = ( 1 - \frac{2}{E} |K|^2 ) T -
\label{e:tbot}
\end{equation}
\[ - \frac{1}{E} \{ ( \overline{K}
\frac{\partial F^j}{\partial z} + K \frac{\partial F^j}{\partial
\overline{z}} ) Z_j + ( \overline{K} \frac{\partial
\overline{F}^j}{\partial z} + K \frac{\partial
\overline{F}^j}{\partial \overline{z}} ) \overline{Z}_j \} .
\] \label{l:8}
\end{lemma}
\noindent {\em Proof}. The characteristic direction decomposes as
$T = \Psi_* T^T + T^\bot$, where $T^T = \lambda \partial /\partial
z + \overline{\lambda}
\partial /\partial \overline{z}$, for some $\lambda \in C^\infty
(U)$. Taking the inner product with $\Psi_* \partial /\partial
\overline{z}$ yields $\lambda = \overline{K}/E$ hence
(\ref{e:partials}) yields (\ref{e:tbot}). $\square$
\begin{lemma} Let $\Psi : N \to {\mathbb H}_n$ be an isometric
immersion of $(N , g)$ into $({\mathbb H}_n , g_0 )$. Then
\begin{equation}
2 \sum_{j=1}^n \frac{\partial F^j}{\partial z} \frac{\partial
F^j}{\partial \overline{z}} + K^2 = 0, \label{e:46}
\end{equation}
\begin{equation}
\sum_{j=1}^n \left( \left| \frac{\partial F^j}{\partial z}
\right|^2 + \left| \frac{\partial F^j}{\partial \overline{z}}
\right|^2 \right) + |K|^2 \neq 0. \label{e:47}
\end{equation} \label{l:9}
\end{lemma}
\noindent {\em Proof}. A calculation based on (\ref{e:partials})
shows that the Webster metric of $({\mathbb H}_n , \theta_0 )$ is
given (with respect to the frame $\{ \partial /\partial x^j , \;
\partial /\partial y^j , \; \partial /\partial t \}$) by
\[ g_0 : \left( \begin{array}{ccc}
2 (\delta_{jk} + 2 y_j y_k ) & - 4 y_j x_k & - 2 y_j \\
- 4 x_j y_k & 2 (\delta_{jk} + 2 x_j x_k ) & 2 x_j \\
- 2 y_k & 2 x_k & 1 \end{array} \right)  \] hence
\[ g_\theta (\Psi_* \frac{\partial}{\partial z} \, , \, \Psi_*
\frac{\partial}{\partial \overline{z}}) = \Psi^A_z
\Psi^B_{\overline{z}} g_{AB} = |K|^2 + \sum_j (|F^j_z |^2 +
|F^j_{\overline{z}}|^2 ), \]
\[ g_\theta (\Psi_* \frac{\partial}{\partial z} \, , \, \Psi_*
\frac{\partial}{\partial z}) = \Psi^A_z \Psi^B_z g_{AB} = K^2 +
\sum_j F^j_z F^j_{\overline{z}} ,  \] (where $g_{AB} = g_0
(\partial_A , \partial_B )$). Since $\Psi$ is an isometric
immersion
\begin{equation}
g_0 (\Psi_* \frac{\partial}{\partial x} \; ,  \; \Psi_*
\frac{\partial}{\partial y}) = 0, \label{e:48}
\end{equation}
\begin{equation}
g_0 (\Psi_* \frac{\partial}{\partial x} \; , \; \Psi_*
\frac{\partial}{\partial x}) = g_0 (\Psi_*
\frac{\partial}{\partial y} \; , \; \Psi_*
\frac{\partial}{\partial y}), \label{e:49}
\end{equation}
and then (\ref{e:48})-(\ref{e:49}) yield
(\ref{e:46})-(\ref{e:47}), respectively. $\square$ \vskip 0.1in
Note that (again by (\ref{e:partials}))
\[ \Delta \Psi = (\Delta \psi^A ) \partial_A = (\Delta F^j ) Z_j +
(\Delta \overline{F}^j ) \overline{Z}_j + \]
\[ +  \{ \Delta f + 2 \sum_{j=1}^n (\Psi^{j} \Delta \Psi^{j+n} -
\Psi^{j+n} \Delta \Psi^{j}) \} T \] and (by Lemma \ref{l:8})
\[ i E \, J T^\bot =  (\overline{K} F^j_z + K F^j_{\overline{z}})
Z_j - (\overline{K} \; \overline{F}^j_z + K
\overline{F}^j_{\overline{z}}) \overline{Z}_j \] hence the
minimality condition (\ref{e:delta}) becomes
\begin{equation}
\Delta F^j = - \frac{2i}{E} (\overline{K} F^j_z + K
F^j_{\overline{z}}), \;\;\; 1 \leq j \leq n, \label{e:delf}
\end{equation}
and $\Delta f = \frac{i}{2} \sum_j (\overline{F}^j \Delta F^j -
F^j \Delta \overline{F}^j )$ or (by (\ref{e:delf}))
\begin{equation}
\Delta f = \frac{1}{E} \{ \overline{K} (|F|^2 )_z + K  (|F|^2
)_{\overline{z}} \} .
\label{e:delfp}
\end{equation}
Let $N$ be a Riemann surface. An immersion $\Psi : N \to {\mathbb
H}_n$ is {\em conformal} if (\ref{e:48})-(\ref{e:49}) hold, for
any local complex coordinate system $(U , z = x + i y)$ on $N$.
Moreover (\ref{e:delf})-(\ref{e:delfp}) lead to the following
definition. A {\em minimal surface} in ${\mathbb H}_n$ is a
Riemann surface $N$ together with a conformal immersion $\Psi : N
\to {\mathbb H}_n$ such that
\begin{equation}
F^j_{z\overline{z}} + i (\overline{K} F^j_z + K
F^j_{\overline{z}}) = 0, \;\;\; 1 \leq j \leq n, \label{e:delf1}
\end{equation}
\begin{equation}
f_{z\overline{z}} - \frac{1}{2} \{ \overline{K} (|F|^2 )_z + K
(|F|^2 )_{\overline{z}} \} = 0. \label{e:delfp1}
\end{equation}
Here $|F|^2 = \sum_j F^j \overline{F}^j$. We may state the
following
\begin{theorem} Let $\Omega \subset {\mathbb C}$ be a simply
connected domain and $\Psi : \Omega \to {\mathbb H}_n$ a minimal
surface such that $J T^\bot = 0$ {\rm (}e.g. $\Psi (\Omega )$ is
tangent to the characteristic direction of $d \theta_0${\rm )}.
Let us set $\Phi =
\partial \Psi /\partial z$. Then $\Phi$ is holomorphic and {\rm
(\ref{e:46})-(\ref{e:47})} hold in $\Omega$. Viceversa, let $\Phi
: \Omega \to {\mathbb C}^{2n+1}$ be a holomorphic map and let us
set
\begin{equation}
\label{e:psi} \Psi^A (z) = {\rm Re} \int_o^z \Phi^j (\zeta ) d
\zeta , \;\;\; A \in \{ 0, 1, \cdots , 2n \} ,
\end{equation}
for any $z \in \Omega$, where $o \in \Omega$ is a fixed base
point. Let $K : \Omega \to {\mathbb C}$ be given by \[ K = \Phi^0
- 2 \sum_{j=1}^n \{ \Phi^j \; {\rm Re} \int_o^z \Phi^{j+n}(\zeta )
d \zeta + \Phi^{j+n} \; {\rm Re} \int_o^z \Phi^j (\zeta ) d \zeta
\} .
\]
If the following identities hold in $\Omega$
\begin{equation} 2 \sum_{j=1}^n \{ |\Phi^j |^2 - |\Phi^{j+n}|^2 +
i(\Phi^{j+n} \overline{\Phi}^j + \Phi^j \overline{\Phi}^{j+n} ) \}
+ K^2 = 0, \label{e:new55}
\end{equation}
\begin{equation}
2 \sum_{j=1}^n (|\Phi^j |^2 + |\Phi^{j+n}|^2 ) + |K|^2 \neq 0,
\label{e:new56}
\end{equation}
\begin{equation}
\overline{K} (\Phi^j + i \Phi^{j+n}) + K (\overline{\Phi}^j + i
\overline{\Phi}^{j+n}) = 0, \;\;\; 1 \leq j \leq n, \label{e:54}
\end{equation}
then $\Psi : \Omega \to {\mathbb H}_n$ is a minimal immersion such
that $J T^\bot = 0$. \label{t:7}
\end{theorem}
\noindent Compare to Theorem 8.1 in \cite{kn:DW}, p. 58. {\em
Proof of Theorem} \ref{t:7}. (\ref{e:46})-(\ref{e:47}) follow from
Lemma \ref{l:9}. Next $J T^\bot = 0$ and
(\ref{e:delf})-(\ref{e:delfp}) yield $\partial \Phi /\partial
\overline{z} = 0$ in $\Omega$.
\par
Viceversa, given a holomorphic map $\Phi : \Omega \to {\mathbb
C}^{2n+1}$ the function $\Psi^A$ given by (\ref{e:psi}) is well
defined (by the classical theorem of Cauchy the integral doesn't
depend upon the choice of path from $o$ to $z$) and $\partial \Psi
/\partial z = \Phi$ hence (\ref{e:new55})-(\ref{e:new56}) yield
(\ref{e:46})-(\ref{e:47}) so that (\ref{e:48})-(\ref{e:49}) are
satisfied and $g_0 (\Psi_*
\partial /\partial x \, , \, \Psi_* \partial /\partial x) \neq 0$,
i.e. $\Psi$ is a conformal immersion. Finally (\ref{e:54}) may be
written
\[ \overline{K} F^j_z + K F^j_{\overline{z}} = 0, \;\;\; 1 \leq j
\leq n, \] which is equivalent (by Lemma \ref{l:8}) to $J T^\bot =
0$ and (\ref{e:delf1})-(\ref{e:delfp1}) imply minimality.
$\square$

\end{document}